%% file: 1995-106.tex
\documentstyle[12pt,amssymbt]{amsart}
\input{lastart}

\begin{document}

\newcommand{\vse}{\vspace{.2in}}
\renewcommand{\theequation}{\thesection.\arabic{equation}}
\newcommand{\shd}{\hspace{-2.5ex plus  -.2ex}.\hspace{1.4ex plus .1ex}}


\newcommand{\thd}{\hspace{-1.1ex}{\bf.}\hspace{1.1ex}}

\title
[\hskip-1em Approximation by analytic matrix 
functions.~The four block problem]
{Approximation by analytic matrix functions.\\
The four block problem}

\author{Vladimir. V. Peller}
\address{\hskip-\parindent
Vladimir Peller\\
Department of Mathematics\\
Kansas State University\\
Manhattan, Kansas 66506}
\email{peller@@math.ksu.edu}

\author{Sergei R. Treil}
\address{\hskip-\parindent Sergei Treil\\
Department of Mathematics\\
Michigan State University\\
East Lansing, Michigan 48824}

\email{treil@@math.msu.edu}
\subjclass{47B35, 93B36}

\thanks{Both authors are partially supported by NSF grant DMS 9304011.
Research at MSRI is supported in part by NSF grant DMS-9022140.}

\begin{abstract}
We study the problem of finding a superoptimal solution to the four
block problem. Given a bounded block matrix function
$\left(\begin{array}{cc}\Phi_{11}
&\Phi_{12}\\\Phi_{21}&\Phi_{22}\end{array}\right)$ on the unit circle
the four block problem is to minimize the $L^\infty$ norm of
$\left(\begin{array}{cc}
\Phi_{11}-F&\Phi_{12}\\\Phi_{21}&\Phi_{22}\end{array}\right)$ over
$F\in H^\infty$. Such a minimizing $F$ (an optimal solution) is almost
never unique. We consider the problem to find a superoptimal solution
which minimizes not only the supremum of the matrix norms but also the
suprema of all further singular values. We give a natural condition
under which the superoptimal solution is unique.
\end{abstract}

\maketitle

\newcommand{\ov}{\overline}
\newcommand{\W}{{\cal W}}
\renewcommand{\i}{\iota}

\section{\shd Introduction}
\label{s1}
\setcounter{equation}{0}

\

The problem of approximating a given scalar function $\f$ on the unit circle
$\T$ uniformly by functions analytic in the unit disk $\dd$ has been attracting
analysts for a long time (see [Kha], [RSh], [Ne], [AAK1-2], [CJ], [PKh]).
It was shown in [Kha] that for a continuous function $\f$ such a best approximation
is unique while it is not unique in the general case. Later it turned out
that this problem is closely related with Hankel operators. Namely, it was proved
by Nehari [Ne] that 
$$
\dist_{L^\be}(\f,H^\be)=\|H_\f\|,
$$
where $H_\f:\;H^2\to H^2_-\df L^2\ominus H^2$ is the {\it Hankel operator} with symbol
$\f$ defined by
$$
H_\f f=\pp_-\f f,~~~f\in H^2,
$$
(we denote by $\pp_+$ and $\pp_-$ the orthogonal projections onto $H^2$ and 
$H^2_-$). Presently the problem of approximating by analytic functions in $L^\be$
is called {\it Nehari's problem}.

We shall also need the notion of a 
Toeplitz operator. Given $\f\in L^\be$ the {\it Toeplitz operator}
$T_\f:\;H^2\to H^2$ is defined by
$$
T_\f f=\pp_+\f f,~~~f\in H^2.
$$

Adamyan, Arov and Krein [AAK1-2] found many interesting connections between Hankel
operators and Nehari's problem. In particular they found 
a more general condition
under which a best approximation is unique: if the essential norm 
$\|H_\f\|_{\rm e}$ is less than $\|H_\f\|$, then $\f$ has a unique best 
approximation. (Recall that for an operator $T$ on Hilbert space
$$
\|T\|_{\rm e}\df\inf\{\|T-K\|:~K~\mbox{is compact}\}). 
$$
This sufficient uniqueness condition can easily be reformulated in terms of
the function $\f$ itself since
$$
\|H_\f\|_{\rm e}=\dist_{L^\be}(\f,H^\be+C)
$$
(see [AAK1-2]).
They also found a criterion of uniqueness
of a best approximation in terms of the corresponding Hankel operator, and in 
the case of non-uniqueness parametrized all
best approximations (optimal solutions of Nehari's problem), see
[AAK1-2]. However it is 
not very easy to verify whether a function $\f$ in $L^\be$ satisfies 
the criterion.
 
Carleson and Jacobs [CJ] studied smoothness properties of the best approximation
for smooth functions $\f$. They proved that if $\f$ belongs to the H\"{o}lder--Zygmund
class $\L_\a$, $\a>0$, $\a\not\in\Z$, then the best approximation also belongs to
the same class. 

Later in [PKh] more general hereditary properties of the non-linear 
operator of best approximation were studied. For a large class of function spaces $X$
on $\T$ it was proved that if $\f\in X$ and $f$ is the best approximation by analytic 
functions, then $f\in X$. Note also that in [PKh] Nehari's problem was also applied
in prediction theory which led to a new approach to the problem of describing
stationary processes satisfying various regularity conditions in terms of their
spectral densities.

A new wave of interest in Nehari's problem was caused by the development of 
$H^\be$ control theory where Nehari's problem plays a central role (see [Fr]).
Moreover for the needs of $H^\be$ control theory it is important to consider
Nehari's problem for matrix-valued functions: given an $n\times m$ matrix function $\Phi$
on $\T$ the problem is to approximate $\Phi$ by bounded analytic matrix functions
$Q$ in the norm
$$
\|\Phi-Q\|_\be\df\ess\sup_{\z\in\T}\|\Phi(\z)-F(\z)\|,
$$
where $\|\cdot\|$ on the right-hand side is the norm of the matrix as an operator from $\C^m$ to $\C^n$.
However in contrast with the scalar case we have uniqueness of a best approximation
only in exceptional cases. Indeed, let
$$
\Phi=\left(\begin{array}{cc}\bar z&0\\0&\frac{1}{2}\bar z\end{array}\right).
$$
Clearly $\dist_{L^\be}(\bar z,H^\be)=1$ and so $\|\Phi-Q\|_\be\ge1$ for any
$Q\in H^\be$. However it is easy to see that any function of the form
$\left(\begin{array}{cc}0&0\\0&q\end{array}\right)$ with $q\in H^\be$,
$\|q\|_\be\le\frac{1}{2}$,is a best approximation. Intuitively, however, is clear
that the ``very best'' approximation is the zero matrix function $\0$.

In [Y] Young suggested imposing the following additional assumptions on approximating
functions. Let $\O_0$ be the set of best approximations:
$$
\O_0=\{Q\in H^\be:~Q~\mbox{minimizes}~~\ess\sup_{\z\in\T}\|\Phi(\z)-Q(\z)\|\}.
$$
Define inductively the sets $\O_j$ as follows
$$
\O_j=\{Q\in \O_{j-1}:
~Q~\mbox{minimizes}~~\ess\sup_{\z\in\T}s_j(\Phi(\z)-Q(\z))\}
$$
(for a matrix (or an operator) $A$ the $j$th {\it singular value} $s_j(A)$, $j\ge0$, is the distance from
$A$ to the set of matrices (operators) of rank at most $j$, $s_0(A)\df\|A\|$).
Elements of $\O_{\min\{m,n\}-1}$ are called {\it superoptimal approximations} 
of $\Phi$ (or superoptimal solutions of Nehari's problem). Put
$$
t_j\df\ess\sup_{\z\in\T}s_j(\Phi(\z)-Q(\z)),~~~Q\in\O_j.
$$
The numbers $t_j$, $0\le j\le\min\{m,n\}-1$, are called the
{\it superoptimal singular values}. As in the scalar case, $t_0=\|H_\Phi\|$
(the Hankel operator $H_\Phi:H^2(\C^n)\to H^2(\C^m)$ is defined in the same way 
as in the scalar case).

Note that $Q$ is a superoptimal solution of Nehari's problem if and only if
it lexicographically minimizes the sequence $\{s_j^\be(\Phi-Q)\}_{j\ge0}$,
where for a matrix function $F$ on $\T$
$$
s_j^\be(F)\df\ess\sup_{\z\in\T}s_j(F(\z)).
$$
 
It was proved in [PY1] that for $\Phi\in H^\be+C$ there exists a unique 
superoptimal approximation $Q$. The method of the proof in [PY1] 
is based on certain
special factorizations of matrix functions (thematic factorizations) and it is
constructive. Later in [T] another method was suggested to establish 
uniqueness
in the $H^\be+C$ case which is based on weighted Nehari's problem.

However in the case when the Hankel operator $H_\Phi$ is non-compact, 
there was no
analog of the Adamyan--Arov--Krein sufficient  condition for uniqueness
in the case of matrix functions.

Nehari's problem is a special case of the so-called 
{\it four block} problem which is one of the most important problems
in control theory. Let $\Phi$ be a block matrix function of the form 
$$
\Phi=\left(
\begin{array}{c|c} \Phi_{11} & \Phi_{12} \\ \hline 
                   \Phi_{21}& \Phi_{22} 
\end{array}
\right).
$$
Here $\Phi$ has size $m\times n$, $\Phi_{11}$ has size $m_1\times n_1$, and
$\Phi_2$ has size $m_2\times n_2$. The four block problem is to minimize
\begin{equation}
\label{1.1}
\left\|\left(
\begin{array}{c|c} \Phi_{11}-Q & \Phi_{12} \\ \hline 
                   \Phi_{21}& \Phi_{22} 
\end{array}
\right)\right\|_\be,
\end{equation}
over bounded analytic functions $Q$ of size $m_1\times n_1$. A function 
$Q\in H^\be(M_{m_1n_1})$ is called an {\it optimal solution} of the four block
problem if it minimizes the norm (1.1).

The four block problem arises naturally when one considers the following (model-matching)
problem in $H^\be$ control. Let $F$, $G_1$ and $G_2$ be matrix 
functions of class $H^\be$. The problem is to minimize
\begin{equation}
\label{1.2}
 \|F-G_1QG_2\|_\be
\end{equation}
over $Q\in H^\be$
(the sizes of the matrix functions in (\ref{1.2}) are such that (1.2) is 
meaningful). Many problems in $H^\be$ control reduce to 
the model-matching problem. 

Engineers usually consider the case of continuous (on  $\T$), or even rational
functions $G_1$ and $G_2$ and assume that these functions
have constant rank on the boundary. Under this assumption 
the model-matching problem reduces
to the four block problem (\ref{1.1}), while it reduces to Nehari's problem
only if the
matrices $G_{2}$, $G_{1}^{*}$ have maximal column rank (the rank equals
the number of columns).
 The assumption on the maximal column rank does not hold for many
interesting applied problems, so engineers have to consider the
four block problem as well.

In the most general case  the model matching problem (\ref{1.2}) 
reduces to the four block problem
under the assumption
that the outer parts of functions $G_1$ and $G_2^*$ are right invertible in $L^\be$ (for continuous functions this is equivalent to 
fact that they have constant
rank on $\T$). The model mathcing problem reduces to Nehari's
problem if the outer parts of functions $G_1$ and $G_2^*$ are  invertible in 
$L^\be$ (which for continuous functions is equivalent to the above maximal column rank assumption).

By analogy with Hankel operators we define the four block operator 
 $\G_\Phi:H^2(\C^{n_1})\oplus L^2(\C^{n_2})\to 
H^2_-(\C^{m_1})\oplus L^2(\C^{m_2})$
by
$$
\G_\Phi\left(\begin{array}{c}f_1\\f_2\end{array}\right)=
\pp^-\Phi\left(\begin{array}{c}f_1\\f_2\end{array}\right),
$$
where $\pp^-$ is the orthogonal projection from $L^2(\C^{m_1})\oplus L^2(\C^{m_2})$
onto  $H^2_-(\C^{m_1})\oplus L^2(\C^{m_2})$. As in the case of 
Hankel operators the {\it infimum} in (\ref{1.1}) is equal to $\|\G_\Phi\|$ 
(see [FT]). The matrix function $\Phi$ is called a {\it symbol} of the four
block operator $\G_\Phi$ (a four block operator has many different symbols). 

As in the case of Nehari's problem we can define the sets $\O_j$:
$$
\O_0=\left\{Q\in H^\be:~Q~\mbox{minimizes}~ 
\left\|\left(
\begin{array}{c|c} \Phi_{11}-F & \Phi_{12} \\ \hline 
                   \Phi_{21}& \Phi_{22} 
\end{array}
\right)\right\|_\be\right\},
$$
$$
\O_j=\left\{Q\in \O_{j-1}:~Q~\mbox{minimizes}~~ 
\ess\sup_{\z\in\T}s_j\left(\left(
\begin{array}{c|c} \Phi_{11}-F & \Phi_{12} \\ \hline 
                   \Phi_{21}& \Phi_{22} 
\end{array}
\right)(\z)\right)\right\}.
$$
$Q$ is called a {\it superoptimal solution} of the four block problem (\ref{1.1}) if
 $Q\in\O_{\min\{m_1,n_1\}-1}$. We define the superoptimal
singular values of the four block problem (\ref{1.1}) by
$$
t_j=\ess\sup_{\z\in\T}s_j\left(\left(
\begin{array}{c|c} \Phi_{11}-Q & \Phi_{12} \\ \hline 
                   \Phi_{21}& \Phi_{22} 
\end{array}
\right)(\z)\right),~~~~F\in\O_j.
$$
Clearly, $t_0=\|\G_\Phi\|$.

We say that a $\Phi$ is a {\it superoptimal symbol} of a four block operator
$\G$ if $\G=\G_\Phi$ and the zero function is a superoptimal solution
of the corresponding four block problem. 

Using a simple compactness argument one can prove easily that a superoptimal
solution always exists.
However we cannot expect a sufficient condition for uniqueness of the superoptimal 
approximation which would be similar to the one found in [PY1] in the case of Nehari's
problem. Indeed it can easily be proved that a four block operator cannot be compact
unless $\Phi_{12},~\Phi_{21}$, and $\Phi_{22}$ are identically equal to zero in
which case the four block problem is equivalent to Nehari's problem.

The main result of the paper (Theorem \ref{t2.3}) is a sufficient condition
for the four block problem to have a unique superoptimal solution. This result
can be considered as an analog of the Adamyan--Arov--Krein theorem mentioned 
above which deals with Nehari's problem in the scalar case. 
Note that Theorem \ref{t2.3} also gives us
a new result for Nehari's problem in the case when the corresponding 
Hankel operator is non-compact. 

The proof is constructive. We give an algorithm to find the unique superoptimal
solution. The algorithm is similar to the one given in [PY1], it reduces the
problem to the case of matrix functions of lower size. However the proof is considerably
more complicated than in the case of Nehari's problem with compact Hankel 
operator. 

In Section 3 we describe briefly the method of factorization and diagonalization.
We construct important matrix functions $V$ and $W$ which can be considered
as analogs of the thematic functions defined in [PY1]. 

In Section 4 we use the construction of Section 3 to parametrize all optimal
solutions of the four block problem. This allows us to reduce the problem  
of finding superoptimal solutions to the case of matrix functions of
lower size.

Section 5 is devoted to the proof of the fact that the matrix functions
$V_c$ and $W_c$ which are submatrices of $V$ and $W$ are left invertible in
$H^\be$. This is one of the principal points in the proof of the main result.

In Section 6 we prove another crucial fact for the proof of the main result. 
Namely, we show that if $\Phi$ satisfies the hypotheses of Theorem \ref{t2.3}, then
the lower order four block problem, obtained as a result of parametrization
in Section 4, also satisfies the hypotheses of Theorem \ref{t2.3}. This
makes it possible to continue the process and complete the proof of the main 
result.

In Section 7 we study superoptimal symbols of four block operators. We
obtain certain special factorizations of such symbols (thematic factorizations),
and define the indices of such factorizations. In the case of Nehari's problem
with compact Hankel operators such factorizations were found in [PY1].

To prove the invariance of indices we introduce in Section 8 the notion of
a superoptimal weight for the four block operator. This is an analog of
the notion introduced in [T] in the case of compact Hankel operators.

In Section [9] we use superoptimal weights to prove that the sums of the indices
in a thematic factorization which correspond to equal superoptimal singular
values do not depend on the choice of factorization. In the case of
compact Hankel operators this invariance property was proved in [PY2].
Note that as in the case of Nehari's problem for $H^\be+C$ functions
if there are equal superoptimal singular values, the indices can depend
on the choice of thematic factorization (see [PY2]).

The last section is devoted to inequalities between the superoptimal singular
values and the singular values of the four block operator. We obtain an
inequality which is new even in the case of Nehari's problem with compact 
Hankel operator. It 
is stronger than the one obtained in [PY2].

Note that in [PY1] hereditary properties of the non-linear operator of
superoptimal approximation were studied. It was shown there that
for a large class of function spaces $X$ the inclusion $\Phi\in X$ implies
that the superoptimal approximant to $\Phi$ also belongs to $X$. It would
be interesting to find analogs of such results in the case of the four block
problem. In particular we do not know whether the superoptimal solution
of the four block problem (\ref{1.1}) must belong to a H\"{o}lder class
$\L_\a$ if $\Phi\in\L_\a$. 

Throughout this paper we shall denote by $M_{m,n}$ the space of $m\times n$
matrices. We shall use the notation $L^\be(M_{m,n})$ and $H^\be(M_{m,n})$
for the spaces of bounded and bounded analytic functions which take
values in $M_{m,n}$. Sometimes if it does not lead to a confusion, we shall simply 
write $L^\be$ and $H^\be$ instead of $L^\be(M_{m,n})$ and $H^\be(M_{m,n})$.

A matrix function $\Theta\in H^\be$ is called {\it inner} if $\Theta(\z)$ 
is isometric
for almost all $\z\in\T$. A matrix function $F\in H^\be (M_{m,n})$
is called {\it outer} if $FH^2(\C^n)$ is dense in $H^2(\C^m)$. A function
$F\in H^\be (M_{m,n})$ is called {\it co-outer} if the transposed function
$F^t\in H^\be(M_{n,m})$ is outer.

\

\

\section{\shd The main result}
\label{s2}
\setcounter{equation}{0}

\

In this section we state the main result of the paper as well as important 
corollaries. Let $\Phi$ be a matrix function of the form
\begin{equation}
\label{2.1}
\Phi=\left(
\begin{array}{c|c} \Phi_{11} & \Phi_{12} \\ \hline 
                   \Phi_{21}& \Phi_{22} 
\end{array}
\right),
\end{equation}
where $\Phi$ has size $m\times n$, $\Phi_{11}$ has size $m_1\times n_1$, and
$\Phi_{22}$ has size $m_2\times n_2$. Recall that $\G_\Phi$ is the four block operator
defined in Section 1 and $\{t_j\}$ is the sequence of superoptimal singular
values. 





The following theorem is the main result of the paper.

\begin{thm}\thd
\label{t2.3}
Let $\Phi$ be a bounded function of the form {\rm (\ref{2.1})}. Suppose that
$\|\G_\Phi\|_{\rm e}$ is less than the smallest nonzero superoptimal singular
value. 
Then there exists a unique superoptimal
solution $Q$ of the four block problem for $\Phi$. The singular values 
$$
s_j\left(\left(
\begin{array}{c|c} \Phi_{11}-Q & \Phi_{12} \\ \hline 
                   \Phi_{21}& \Phi_{22} 
\end{array}
\right)(\z)\right),~~~~1\le j\le d-1,
$$
are constant on $\T$.
\end{thm}


The following partial case of Theorem \ref{t2.3} improves the result of
[PY1] on  uniqueness of superoptimal solutions of Nehari's problem.  

\begin{thm}\thd
\label{t2.4}
Let $\Phi$ be bounded matrix function on $\T$. 
Suppose that $\|H_\Phi\|_{\rm e}$ is less that the smallest nonzero
superoptimal singular value of Nehari's problem. Then 
$\Phi$ has a unique superoptimal approximation $Q$ by bounded analytic matrix
functions. The singular values $s_j(\Phi(\z)-Q(\z))$ are constant on $\T$.
\end{thm}

\

\

\section{\shd Diagonalization}   
\label{s3}
\setcounter{equation}{0}

\

In this section we start with a maximizing vector of the four block operator 
and we construct a certain special unitary-valued matrix. This allows us to
achieve a diagonalization. Later using this diagonalization we shall reduce
the problem to the case of a matrix function of a lower size.

\begin{lem}\thd
\label{t3.1}
Let ${\bf v}$ be an $n\times1$ inner matrix function. Then there exists a co-outer
function $V_c\in H^\be(M_{n,n-1})$ such that the matrix function
$$
{\bf V}\df\left(\begin{array}{cc}{\bf v}&\overline V_c\end{array}\right)
$$
is unitary-valued on $\T$.
\end{lem}

In [PY1] a stronger result was obtained. It was shown that 
all minors of ${\bf V}$ on the first column are in $H^\be$. This property of 
analyticity of minors was essential for the proof of the uniqeness of
a superoptimal solution of Nehari's problem which was given in [PY1].
Earlier the existence of a co-outer $V_c$ satisfying the requirement
of Lemma \ref{t3.1} was proved in [Va], however the property of analyticity
of minors was not noticed in [Va]. It also can be shown that if we 
$V_c^{(1)}$ and $V_c^{(2)}$ are $n\times(n-1)$ co-outer functions
satisfying the requirements of Lemma \ref{t3.1}, then there exists a
constant unitary matrix $U$ such that $V_c^{(1)}=V_c^{(2)}U$ (see [Va], [PY1]).

To start the procedure we need a maximizing vector of the four block operator
$\G_\Phi$, i.e., a nonzero vector $f\in H^2(\C^{n_1})\oplus L^2(\C^{n_2})$ 
such that
$\|\G_\Phi f\|=\|\G_\Phi\|\cdot\|f\|$. 
If $\Phi$ satisfies the hypotheses of Theorem
\ref{t2.3} and $\G_\Phi\neq\0$, then 
$\|\G_\Phi\|=t_0>\|\G_\Phi\|_{\rm e}$
and so a maximizing vector for $\G_\Phi$ exists.

The following fact is well-known in the case of Hankel operators (see [AAK3]).
The proof of it in the case of four block operators is similar.

\begin{lem}\thd
\label{t3.3}
Let $\Phi$ be a matrix function on $\T$ of the form {\rm (\ref{2.1})} and
such that $\|\Phi\|_\be=\|\G_\Phi\|$.
Suppose that $f$ is a maximizing vector for 
$\G_\Phi$ and put  $g=t_0^{-1}\G_\Phi f$. Then 
$\G_\Phi f=\Phi f$ and $\|g(\z)\|_{\C^m}=\|f(\z)\|_{\C^n}$ a.e. on $\T$. 
Furthermore  $\|\Phi(\z)\|=\|\Phi\|_\be$ a.e. on $\T$.
\end{lem}

\pf We have 
$$
\|\G_\Phi f\|_2 =\|{\Bbb P}^-\Phi f\|_2 \le \|\Phi f\|_2\le
\|\Phi\|_\infty\|f\|_2 =
\|\G_\Phi\|\cdot\|f\|_2 = \|\G_\Phi f\|_2.
$$
It follows that all inequalities in this chain are, in
fact, equalities. The fact that $\|\pp^-\Phi\|_2=\|\Phi\|_2$ certainly means
that $\Phi f\in H^2_-(\C^{m_1})\oplus L^2(\C^{m_2})$ and so $\G_\Phi f=\Phi f$.
The equality $\|\Phi f\|_2=\|\Phi\|_\infty\|f\|_2$ implies that
$\|g(\z)\|_{\C^m}=\|f(\z)\|_{\C^n}$ a.e. on $\T$, which in turn implies that 
$\|\Phi(\z)\|=\|\Phi\|_\be$ for almost all $\z\in\T$. $\bl$

\begin{lem}\thd
\label{t3.2}
Let $\Phi$ be a matrix function of the form {\rm (\ref{2.1})}
such that  $\|\G_\Phi\|_{\rm e}<\|\G_\Phi\|$.
Suppose that $f=f_1\oplus f_2$ is a maximizing vector for $\G_\Phi$ and 
 $t_0g=g_1\oplus g_2=\G_\Phi f$, 
where $f_1\in H^2(\C^{n_1})$, $f_2\in L^2(\C^{n_2})$,
$g_1\in H^2_-(\C^{m_1})$, and $g_2\in L^2(\C^{m_2})$. Then
$$
\|f_1(\z)\|^2_{\C^{n_1}}\ge
\frac{\|\G_\Phi\|^2-\|\G_\Phi\|^2_{\rm e}}{\|\G_\Phi\|^2}\|f(\z)\|^2_{\C^n},~~~
\|g_1(\z)\|^2_{\C^{m_1}}\ge
\frac{\|\G_\Phi\|^2-\|\G_\Phi\|^2_{\rm e}}{\|\G_\Phi\|^2}\|g(\z)\|^2_{\C^m}
$$
almost everywhere on $\T$.
\end{lem}

\pf By subtracting an optimal solution, we can assume that 
$\|\Phi\|_\be=\|\G_\Phi\|$. By Lemma \ref{t3.3}, $t_0g(\z)=\Phi(\z)f(\z)$
and $\|g(\z)\|_{\C^m}=\|f(\z)\|_{\C^n}$ a.e. on $\T$. Therefore $f(\z)$ is
a maximizing vector for $\Phi(\z)$ and $g(\z)$ is a maximizing vector for
$\Phi^*(\z)$. We have 
$g_2=\left(\begin{array}{cc}\Phi_{21}&\Phi_{22}\end{array}\right)f$. It is well
known that
$$
\left\|\G_{\left(\begin{array}{c|c}\0&\0\\\hline\Phi_{21}&\Phi_{22}
\end{array}\right)}\right\|_{\rm e}
=\left\|\G_{\left(\begin{array}{c|c}\0&\0\\\hline\Phi_{21}&\Phi_{22}
\end{array}\right)}\right\|
=\left\|\left(\begin{array}{cc}\Phi_{21}&\Phi_{22}\end{array}\right)\right\|_\be
$$
and so
$$
\left\|\left(\begin{array}{cc}\Phi_{21}&\Phi_{22}\end{array}\right)\right\|_\be\le
\|\G_\Phi\|_{\rm e}.
$$
Hence
$$
\|\G_\Phi\|\cdot\|g_2(\z)\|_{\C^{m_2}}\le
\left\|\left(\begin{array}{cc}\Phi_{21}&\Phi_{22}\end{array}\right)\right\|_\be
\|f(\z)\|_{\C^n}\le\|\G_\Phi\|_{\rm e}\|f(\z)\|_{\C^n}.
$$
On the other hand
$$
\|g(\z)\|^2_{\C^m}=\|g_1(\z)\|^2_{\C^{m_1}}+\|g_2(\z)\|^2_{\C^{m_2}}
=\|f(\z)\|_{\C^n}^2.
$$
Therefore
$$
\|g_1\|_{\C^{m_1}}^2\ge
\frac{\|\G_\Phi\|^2-\|\G_\Phi\|^2_{\rm e}}{\|\G_\Phi\|^2}\|g(\z)\|^2_{\C^m}.
$$

To prove the inequality for $f_1(\z)$ we can use the same argument
since 
$t_0f(\z)=\Phi^*(\z)g(\z)$  and
$$
\left\|\left(\begin{array}{cc}\Phi^*_{12}&\Phi^*_{22}\end{array}\right)\right\|=
\left\|\left(\begin{array}{c}\Phi_{12}\\\Phi_{22}\end{array}\right)\right\|_\be=
\left\|\G_{\left(\begin{array}{c|c}\0&\Phi_{12}\\\hline\0&\Phi_{22}\end{array}
\right)}\right\|_{\rm e}
=\left\|\G_{\left(\begin{array}{c|c}\0&\Phi_{12}\\\hline\0&\Phi_{22}\end{array}
\right)}\right\|.~\bl
$$

\begin{cor}\thd
\label{t3.4}
Let $\Phi$ be a matrix function of the form {\rm (\ref{2.1})} such that
 $\|\G_\Phi\|_{\rm e}<\|\G_\Phi\|$.
Suppose that $\min\{m_1,n_1\}=1$. Then there exists a unique optimal solution of
the four block problem.
\end{cor}

\pf It is sufficient to consider the case $n_1=1$. To obtain the result in 
the case
$m_1=1$ we can pass to the transpose of $\Phi$. Let $f=f_1\oplus 
f_2$ be a maximizing
vector of $\G_\Phi$ and let $Q\in H^\be(M_{m_11})$ be
an optimal solution of the four block problem, i.e.,
$$
\left\|\left(\begin{array}{cc}\Phi_{11}-Q&\Phi_{12}\\\Phi_{21}&\Phi_{22}
\end{array}\right)\right\|=\|\G_\Phi\|.
$$
It follows from Lemma \ref{t3.2} that $f_1$ is a nonzero scalar function in $H^2$.
By Lemma \ref{t3.3}
$$
\left(\begin{array}{cc}\Phi_{11}-Q&\Phi_{12}\\\Phi_{21}&\Phi_{22}
\end{array}\right)\left(\begin{array}{c}f_1\\f_2\end{array}\right)=\G_\Phi f.
$$
Therefore $Qf_1$ is uniquely determined by $\G_\Phi$ and $f$
($Q$ is a column matrix) and since $f_1\neq\0$, it follows that $Q$ is
uniquely determined by $\Phi$. $\bl$

Now, following ideas of [PY1] we construct diagonalizing matrix 
functions
$V$ and $W$ in the following way. 

Take a maximizing vector $f=f_1\oplus f_2$,
$f_1\in H^2(\C^{n_1})$, $f_2\in L^2(\C^{n_2})$. Let $h$ be a scalar outer function such that
$|h(\z)|^2=\|f(\z)\|_{\C^n}^2$, $\z\in\T$. Such a function $h$ always exists since $f_1\in 
H^2$ and
$f_1\ne\0$. Denote by $\vt$  a greatest common inner divisor of the
entries of
$f_1$ (it may happen that $\vt={\bf 1}$).  Define the vector $v=v_1\oplus 
v_2\in 
H^\infty(\C^{n_1})\oplus L^\infty(\C^{n_2})$ by $v=\bar\vt f/h$. 
Clearly, $\|v(z)\|=1$ a.e. on $\T$. The vector $v$ will be the first 
column of the matrix $V$. 

We can represent the column function $v_1$ as $v_1=v_{(o)}v_{(i)}$
where $v_{(o)}$ is a scalar outer function such that 
$|v_{(o)}(\z)|=\|v_1(\z)\|_{\C^{n_1}}$
a.e. on $\T$
and $v_{(i)}$ is an inner column function (the inner part of $v_1$), i.e.,
$\|v_{(i)}(\z)\|_{\C^{n_1}}=1$ a.e. on $\T$.
 Applying Lemma \ref{t3.1} 
to ${\bf v}=v_{(i)}$,
we obtain an inner and co-outer matrix $V_c$ such that the matrix
$\left(\begin{array}{cc}v_{(i)}&\overline V_c\end{array}\right)$ 
is unitary-valued. Note that the vector-function 
$v_{(i)}$
is pointwise orthogonal to any column of $\overline V_c$ a.e. on $\T$, and so 
the same is true for $v_1$. So the matrix function
\begin{equation}
\label{3.1}
\left(\begin{array}{cc}
v_1 & \overline V_c \\ 
v_2 &  \0
\end{array}
\right)
\end{equation}
is isometric almost everywhere on $\T$. 

It is easy to see that we can 
complete this matrix function by adding $m_2$ measurable column functions 
to obtain a unitary-valued function. Indeed it is sufficient to complete
the matrix function to a square matrix function whose columns are
pointwise linear independent and then apply the Gram-Schmidt orthogonalization
process to the columns. To this end
we can approximate our matrix
function uniformly by step functions which take isometric values. 
Clearly, we can
find a unitary completion for each step function. It is easy to see that
if the distance from a step function to our initial function is sufficiently
small, then the columns of our initial function and the columns we added to
the step function are linearly independent.

Let $V$ be a unitary-valued completion of the matrix (\ref{3.1}). Then
$V$ has the form
$$
V=\left(
\begin{array}{cc|c}
	v_1 & \overline V_c & \star  \\
	\hline
	v_2 & \0 & \star  
\end{array}
\right)
$$

Let us now construct a unitary-valued matrix $W$ in a similar way. Let 
 $t_0g=\G_\Phi f$. Then by Lemma \ref{t3.3}, 
$\|g(\z)\|=\|f(\z)\|=|h(\z)|$,
$\z\in\T$. Let $\t$ be a greatest common inner divisor of all entries of 
$\bar z\overline g_1$ (recall that $g=g_1\oplus g_2$, $g_1\in 
H^2_-(\C^{m_1})$, $g_2\in L^2(\C^{m_2})$). Define the column function 
$w=w_1\oplus w_2\in H^\be(\C^{m_1})\oplus L^\be(\C^{m_2})$ by 
$w\df\bar z\bar\t\bar g/h$. 
By analogy with (\ref{3.1}) we can find a co-outer
matrix function $W_c$ such that the matrix function $\left(
\begin{array}{cc}
	w_1 & \overline W_c \\
	w_2 &  \0
\end{array}\right)$
takes isometric values on $\T$. We can complete this matrix function to
a unitary-valued matrix function and define $W$ to be its transpose: 
$$
W^t=\left(
\begin{array}{cc|c}
	w_1 & \overline W_c & \star  \\
	\hline
	w_2 & \0 & \star
\end{array}
\right).
$$

To prove Theorem \ref{t2.3} we shall proceed as follows. Let $Q_0$ be
an optimal solution of the four block problem. We shall
prove in the next section that
the matrix function $\left(
\begin{array}{cc} \Phi_{11}-Q_0 & \Phi_{12} \\
                   \Phi_{21}& \Phi_{22} 
\end{array}
\right)$
admits a representation
$$
\left(
\begin{array}{c|c} \Phi_{11}-Q_0 & \Phi_{12} \\ \hline 
                   \Phi_{21}& \Phi_{22} 
\end{array}
\right)=
 W^*
\left(\begin{array}{cc|c} t_0u_0 & \0 & \0\\ \0 & \Phi^{(1)}_{11} &
\Phi^{(1)}_{12}\\[.3pc] \hline & & \\[-.8pc] \0 & \Phi^{(1)}_{21} & \Phi^{(1)}_{22}
\end{array}\right) V^*,
$$
where $u_0=\bar z\bar\vt\bar\t\bar h/h$ and $\Phi^{(1)}_{11}$ is a matrix 
function of size $(m_1-1)\times(n_1-1)$. We shall also prove in Section 4 
that if $Q$ is another
optimal solution, then
$$
\left(
\begin{array}{c|c} \Phi_{11}-Q & \Phi_{12} \\ \hline
                   \Phi_{21}& \Phi_{22} 
\end{array}
\right)=
 W^*
\left(\begin{array}{cc|c} t_0u_0 & \0 & \0\\ \0 & \Phi^{(1)}_{11}-Q_1 &
\Phi^{(1)}_{12}\\[.3pc] \hline & & \\[-.8pc] \0 & \Phi^{(1)}_{21} & \Phi^{(1)}_{22}
\end{array}\right) V^*,
$$
where $Q_1\in H^\be(M_{m_1-1,n_1-1})$ and
$$
\left\|\left(\begin{array}{cc}
 \Phi^{(1)}_{11}-Q_1 & \Phi^{(1)}_{12}\\&  \\[-1pc]
\Phi^{(1)}_{21} & \Phi^{(1)}_{22}
\end{array}\right)\right\|_\infty \le t_0.
$$
Since $V$ and $W$ are unitary-valued, it is easy to see that $Q$ is a superoptimal
solution to the four block problem for the matrix function $\Phi$
if and only if $Q_1$ is a superoptimal solution to the four block problem
for the matrix function
$$
\Phi^{(1)}=
\left(
\begin{array}{c|c} \Phi^{(1)}_{11} & \Phi^{(1)}_{12} \\ \hline &  \\[-1pc]
                   \Phi^{(1)}_{21}& \Phi^{(1)}_{22} 
\end{array}
\right).
$$
Moreover, if $t_0,t_1,\cdots,t_{d-1}$ is sequence of superoptimal singular values
of the four block problem for $\Phi$, then $t_1,\cdots,t_{d-1}$ is the sequence of
superoptimal singular values of the four block problem for $\Phi^{(1)}$.
This reduction allows us to diminish the size of the matrix function $\Phi_{11}$.

If $\G_{\Phi^{(1)}}=\0$, we clearly have uniqueness.
To continue this process we have to be able to find a maximizing vector for the
four block operator $\G_{\Phi^{(1)}}$.
We can certainly do that if its
essential norm is still less than 
the smallest nonzero superoptimal singular value. 
In Section 6 we shall prove that
$\|\G_{\Phi^{(1)}}\|_{\rm e}\le\|\G_\Phi\|_{\rm e}$ 
which will allow us to continue
the process and reduce Theorem \ref{t2.3} to Corollary \ref{t3.4}.
  
\

\

\section{\shd Parametrization of optimal solutions}
\label{s4}
\setcounter{equation}{0}

\

In this section we describe the optimal solutions of the four block problem
in case when $\|\G_\Phi\|_{\rm e}<\|\G_\Phi\|$.

\begin{lem}\thd
\label{t4.1}
Let $\Phi$ be a block matrix function of the form {\rm (\ref{2.1})}
such that $\|\G_\Phi\|_{\rm e}<\|\G_\Phi\|$, 
and let $V$ and $W$ be the matrix functions constructed in 
Section 3. Then there exists a unimodular function $u_0$ such that any optimal
solution $Q_0$ of the four block problem satisfies
\begin{equation}
\label{4.1}
\left(
\begin{array}{c|c} \Phi_{11}-Q_0 & \Phi_{12} \\ \hline 
                   \Phi_{21}& \Phi_{22} 
\end{array}
\right)=
 W^*
\left(\begin{array}{cc|c} t_0u_0 & \0 & \0\\ \0 & \Phi^{(1)}_{11} &
\Phi^{(1)}_{12}\\[.3pc] \hline & & \\[-.8pc] \0 & \Phi^{(1)}_{21} & \Phi^{(1)}_{22}
\end{array}\right) V^*,
\end{equation}
where 
$\Phi^{(1)}_{11}$ 
is a matrix function of size $m_1-1\times n_1-1$.
The unimodular function $u_0$ admits a representation 
$u_0=\bar z\bar b\bar h/h$, where $h$ is an outer function in $H^2$ and $b$ 
is a finite Blaschke product. Moreover, the Toeplitz operator $T_{u_0}$ is
Fredholm and $\Range T_{u_0}=H^2$.
\end{lem}

\pf Let $f=f_1\oplus f_2$ be a maximizing vector for $\G_\Phi$ and let
$\G_\Phi f=t_0g=t_0(g_1\oplus g_2)$.
Put $u_0=\bar z\bar\vt\bar\t\bar h/h$ (see the construction of the
matrix functions $V$ and $W$ in Section 3).
By Lemma \ref {t3.3},
$f(\z)$ is a maximizing vector for 
$\left(
\begin{array}{cc} \Phi_{11}(\z)-Q_0(\z) & \Phi_{12}(\z) \\  
                   \Phi_{21}(\z)& \Phi_{22}(\z) 
\end{array}
\right)$
almost everywhere on $\T$ and
$$
\left(
\begin{array}{cc} \Phi_{11}-Q_0 & \Phi_{12} \\ 
                   \Phi_{21}& \Phi_{22} 
\end{array}
\right)f=t_0g.
$$
Therefore $g(\z)$ is a maximizing vector for
$\left(
\begin{array}{cc} \Phi_{11}(\z)-Q_0(\z) & \Phi_{12}(\z) \\  
                   \Phi_{21}(\z)& \Phi_{22}(\z) 
\end{array}
\right)^*$
for almost all $\z\in\T$ and so
\begin{equation}
\label{4.2}
\left(
\begin{array}{cc} \Phi_{11}-Q_0 & \Phi_{12} \\ 
                   \Phi_{21}& \Phi_{22} 
\end{array}
\right)^*g=t_0f.
\end{equation}
Since $f=\vt h\left(\begin{array}{c}v_1\\v_2\end{array}\right)$ and
$\bar z\bar g=\t h\left(\begin{array}{c}w_1\\w_2\end{array}\right)$,
we have
$$
\vt h\left(
\begin{array}{cc} \Phi_{11}-Q_0 & \Phi_{12} \\ 
                   \Phi_{21}& \Phi_{22} 
\end{array}
\right)
\left(\begin{array}{c}v_1\\v_2\end{array}\right)
=t_0\bar z\bar\t\bar h\left(\begin{array}{c}\overline w_1\\\overline w_2\end{array}\right).
$$

It follows from the definition of the matrix functions $V$ and $W$ (see Section 3)
that
$$
\vt h\left(
\begin{array}{cc} \Phi_{11}-Q_0 & \Phi_{12} \\ 
                   \Phi_{21}& \Phi_{22} 
\end{array}
\right)V
\left(\begin{array}{c}{\bf 1}\\\0\\\vdots\\\0\end{array}\right)
=t_0\bar z\bar\t\bar hW^*\left(\begin{array}{c}{\bf 1}\\\0\\\vdots\\\0\end{array}\right).
$$
It is easy to see that the first column of $W\left(
\begin{array}{cc} \Phi_{11}-Q_0 & \Phi_{12} \\ 
                   \Phi_{21}& \Phi_{22} 
\end{array}
\right)V$
has the form 
 $\left(\begin{array}{cccc}t_0u_0&\0&\cdots&\0\end{array}\right)^t$,
where $u_0\df\bar z\bar\vt\bar\t\bar h/h$. Similarly, using (\ref{4.2})
we find that the first row of $W\left(
\begin{array}{cc} \Phi_{11}-Q_0 & \Phi_{12} \\ 
                   \Phi_{21}& \Phi_{22} 
\end{array}
\right)V$ has the form 
$\left(\begin{array}{cccc}t_0u_0&\0&\cdots&\0\end{array}\right)$, which
proves that $\left(
\begin{array}{cc} \Phi_{11}-Q_0 & \Phi_{12} \\ 
                   \Phi_{21}& \Phi_{22} 
\end{array}
\right)$ has the form (\ref{4.1}).

Let us show that the Toeplitz operator $T_{u_0}$ is Fredholm and is onto.
Clearly, $\|H_{u_0}\|=1$, since $\|H_{u_0}h\|_2=\|\pp_-\bar z\bar\vt\bar\t\|_2
=\|\bar z\bar\vt\bar\t\|_2=\|h\|_2$. We claim that 
 $\|H_{u_0}\|_{\rm e}<1$. Indeed, let $f$ be a scalar function in
$H^2$. We have
\bay
\G_\Phi vf&=&\pp_-\left(
\begin{array}{cc} \Phi_{11}-Q_0 & \Phi_{12} \\ 
                   \Phi_{21}& \Phi_{22} 
\end{array}
\right)\left(\begin{array}{c}v_1f\\v_2f\end{array}\right)\nonumber\\[.5pc]
&=&\pp_-W^*
\left(\begin{array}{ccc} t_0u_0 & \0 & \0\\ \0 & \Phi^{(1)}_{11} &
\Phi^{(1)}_{12}\\ [-1pc]&\\ \0 & \Phi^{(1)}_{21} & \Phi^{(1)}_{22}
\end{array}\right) V^*\left(\begin{array}{c}v_1f\\v_2f\end{array}\right)\nonumber\\
&=&\pp_-W^*
\left(\begin{array}{ccc} t_0u_0 & \0 & \0\\ \0 & \Phi^{(1)}_{11} &
\Phi^{(1)}_{12}\\[-1pc]&\\ \0 & \Phi^{(1)}_{21} & \Phi^{(1)}_{22}
\end{array}\right)\left(\begin{array}{c}f\\\0\\\vdots\\\0\end{array}\right)
=t_0\pp_-\left(\begin{array}{c}\overline w_1u_0f\\\overline w_2u_0f\end{array}\right)\nonumber
\ey
Therefore
$$
\pp_-w^t\G_\Phi vf
=t_0\pp_-w^t\pp_-
\left(\begin{array}{c}\overline w_1u_0f\\\overline w_2u_0f\end{array}\right)
=t_0\pp_-\left(\begin{array}{cc}w_1&w_2\end{array}\right)
\left(\begin{array}{c}\overline w_1u_0f\\\overline w_2u_0f\end{array}\right)
=t_0H_{u_0}f,
$$
whence
$$
t_0\|H_{u_0}\|_{\rm e}\le\|v\|_2\|w\|_2\|\G_\Phi\|_{\rm e}=\|\G_\Phi\|_{\rm e}<t_0,
$$
which implies that $\|H_{u_0}\|_{\rm e}<1$. 

Since $\|H_{u_0}\|_{\rm e}=\dist_{L^\be}(u_0,H^\be+C)$ (see e.g., [S], [Ni]),
it follows that 
 $\|H_{u_0}\|_{\rm e}=\lim_{j\to\be}\dist_{L^\be}(z^ju_0,H^\be)$.
We have $\|H_{u_0}h\|=\|h\|$. So $\dist(u_0,H^\be)=\|H_{u_0}\|=1$.
Therefore there exists a $j\in\Z_+$ such that 
$$
\dist_{L^\be}(z^ju_0,H^\be)=1\quad\hbox{and}\quad
\dist_{L^\be}(z^{j+1}u_0,H^\be)<1.$$ This means that $T_{z^{j+1}u_0}$
is left invertible and $T_{z^ju_0}$ is not left invertible which implies that
$T_{z^{j+1}u_0}$ is invertible (see [Ni]). Clearly,
$T_{z^{j+1}u_0}=T_{u_0}T_{z^{j+1}}$, $T_{z^{j+1}}$ is Fredholm and so is
$T_{u_0}$.

Since $u_0$ has the form 
$u_0=\bar z\bar\vt\bar\t\bar h/h$, where $\vt$ and $\t$ are inner and $h$ is
an outer function in $H^2$, the Toeplitz operator has dense range (see [PKh]) which 
together with the Fredholmness of $T_{u_0}$ implies that $T_{u_0}$ is onto. 

It remains to show that  both $\vt$ and $\t$ are finite Blaschke products. Indeed,
if $\kappa$ is an inner divisor of $\vt\t$, it is easy to see that
$\kappa h\in\Ker T_{u_0}$ and since $T_{u_0}$ is Fredholm, $\Ker T_{u_0}$
is finite dimensional, which implies that both $\vt$ and $\t$ are finite Blaschke
products. $\bl$

\begin{thm}\thd
\label{t4.2}
Let $\Phi$ be a block matrix function of the form {\rm (\ref{2.1})}
such that $\|\G_\Phi\|_{\rm e}<\|\G_\Phi\|$
and let $Q_0$ be an optimal solution of the four block problem. Suppose that
$V,~W,~u_0,~\Phi^{(1)}_{11},~\Phi^{(1)}_{12},~\Phi^{(1)}_{21},~\Phi^{(1)}_{22}$
satisfy {\rm (\ref{4.1})} holds. Let $Q$ be a matrix function of size $m_1\times n_1$.
Then $Q$ is an optimal solution of the four block problem if and only if
there exists $Q_1\in H^\be(M_{m_1-1,n_1-1})$ that satisfies the following conditions:
\begin{equation}
\label{4.3}
\left(
\begin{array}{c|c} \Phi_{11}-Q & \Phi_{12} \\ \hline 
                   \Phi_{21}& \Phi_{22} 
\end{array}
\right)=
 W^*
\left(\begin{array}{cc|c} t_0u_0 & \0 & \0\\ \0 & \Phi^{(1)}_{11}-Q_1 &
\Phi^{(1)}_{12}\\[.3pc] \hline & & \\[-.8pc] \0 & \Phi^{(1)}_{21} & \Phi^{(1)}_{22}
\end{array}\right) V^*,
\end{equation}
\begin{equation}
\label{4.4}
\left\|\left(\begin{array}{cc}
 \Phi^{(1)}_{11}-Q_1 & \Phi^{(1)}_{12}\\[-1pc]&\\
\Phi^{(1)}_{21} & \Phi^{(1)}_{22}
\end{array}\right)\right\|_\infty \le t_0.
\end{equation}
\end{thm}

To prove Theorem \ref{t4.2} we need the following result from [PY1]:

\begin{lem}\thd
\label{t4.3}
Let ${\bf V},{\bf  W}$ be $L^\infty$ matrix functions on $\T$, of types $n\times n$, 
$m \times m$ respectively, which are unitary-valued a.e. and are of the form
$$
{\bf V} = \left(\begin{array}{cc}{\bf v}&\overline V_c\end{array}\right),~~~ 
{\bf W}^t = \left(\begin{array}{cc}{\bf w}&\overline W_c\end{array}\right),
$$ 
where ${\bf v}, V_c, {\bf w}, W_c$ are $H^\infty$ matrix functions, 
${\bf v}$ and ${\bf w}$ are
column functions, and $V_c,~W_c$ are co-outer.  Then
\[
{\bf W} H^\infty (M_{m,n}) {\bf V}  \bigcap \left( \begin{array}{cc} \0 & \0 \\ \0 &
L^\infty (M_{m-1, n-1}) \end{array} \right) = \left( \begin{array}{cc} \0 & 
\0 \\
\0 & H^\infty (M_{m-1, n-1}) \end{array}
\right). 
\]
\end{lem}

{\bf Proof of Theorem \ref{t4.2}.} Let $Q$ be an optimal solution. By Lemma
\ref{t4.1} 
$$W\left(
\begin{array}{c|c} Q_0-Q & \0 \\ \hline 
                   \0& \0
\end{array}
\right)V$$ 
has the form
$$
\left(\begin{array}{cc|c} \0 & \0 & \0\\ \0 & \star &
\star\\ \hline  \0 & \star & \star
\end{array}\right)
$$
(the upper left block is scalar).
On the other hand it is easy to see from the definition of $V$ and $W$ (see
Section 3) that
$$ W\left(
\begin{array}{c|c} Q_0-Q & \0 \\ \hline 
                   \0& \0
\end{array}
\right) V = 
\left(\begin{array}{c|c}
\left(\begin{array}{cc}w_1&\overline W_c\end{array}\right)^t 
(Q_0-Q)
\left(\begin{array}{cc} v_1&\overline V_c
\end{array}\right) & \0\\[.3pc]
\hline & \\[-.8pc] \0 & \0 \end{array}\right).
$$
Therefore
$$
\left(\begin{array}{cc}w_1&\overline W_c\end{array}\right)^t 
(Q_0-Q)
\left(\begin{array}{cc} v_1&\overline V_c
\end{array}\right)
=\left(
\begin{array}{cc}\0&\0\\\0&F\end{array}\right)
$$
for some $F\in L^\be(M_{m_1-1,n_1-1})$ (the upper left corner of the matrix 
function on the right hand side is scalar). Let $v_1=v_{(o)}v_{(i)}$,
$w_1=w_{(o)}w_{(i)}$, where $v_{(o)}$ and $w_{(o)}$ are scalar outer functions, and
$v_{(i)}$ and $w_{(i)}$ are inner column functions.
We have
$$
\left(\begin{array}{cc}w_{(i)}&\overline W_c\end{array}\right)^t 
(Q_0-Q)
\left(\begin{array}{cc} v_{(i)}&\overline V_c
\end{array}\right)
$$
$$
=\left(\begin{array}{cc}w_{(o)}&\0\\\0&I\end{array}\right)
\left(\begin{array}{cc}w_1&\overline W_c\end{array}\right)^t 
(Q_0-Q)
\left(\begin{array}{cc} v_1&\overline V_c
\end{array}\right)\left(\begin{array}{cc}v_{(o)}&\0\\\0&I\end{array}\right)
$$
$$
=\left(\begin{array}{cc}w_{(o)}&\0\\\0&I\end{array}\right)
\left(
\begin{array}{cc}\0&\0\\\0&F\end{array}\right)
\left(\begin{array}{cc}v_{(o)}&\0\\\0&I\end{array}\right)
=\left(
\begin{array}{cc}\0&\0\\\0&F\end{array}\right).
$$
Put ${\bf v}=v_{(i)}$, ${\bf w}=w_{(i)}$.
Clearly, the matrix functions 
${\bf V}=\left(\begin{array}{cc}{\bf v}&\overline V_c\end{array}\right)$ and
 ${\bf W}=\left(\begin{array}{cc}{\bf w}&\overline W_c\end{array}\right)$ 
satisfy the hypotheses of Lemma \ref{t4.3}. Therefore 
 $F\in H^\be(M_{m_1-1,n_1-1})$, 
which proves that $Q$ satisfies (\ref{4.3})
with $Q_1=-F$. Since $Q$ is an optimal solution, (\ref{4.4}) obviously holds.

Conversely, suppose that $Q_1$ is a function in $H^\be(M_{m_1-1,n_1-1})$
satisfying (\ref{4.3}). Then it follows from Lemma \ref{t4.3} that there exists
a function $G\in H^\be(M_{m_1,n_1})$ such that 
$$
\left(\begin{array}{cc}w_1&\overline W_c\end{array}\right)^t 
G\left(\begin{array}{cc} v_1&\overline V_c
\end{array}\right)
=\left(
\begin{array}{cc}\0&\0\\\0&-Q_1\end{array}\right),
$$
which implies that
$$
\left(
\begin{array}{c|c} \Phi_{11}-(G+Q_0) & \Phi_{12} \\ \hline 
                   \Phi_{21}& \Phi_{22} 
\end{array}
\right)=
 W^*
\left(\begin{array}{cc|c} t_0u_0 & \0 & \0\\ \0 & \Phi^{(1)}_{11}-Q_1 &
\Phi^{(1)}_{12}\\[.3pc] \hline & & \\[-.8pc] \0 & \Phi^{(1)}_{21} & \Phi^{(1)}_{22}
\end{array}\right) V^*
$$
and so $Q=G+Q_0\in H^\be(M_{m_1,n_1})$. Clearly, (\ref{4.4}) implies now
that $Q$ is an optimal solution. $\bl$

It is easy to see that Theorem \ref{4.2} reduces the problem of finding a superoptimal solution
for $\Phi$ to the same problem for the matrix function 
$\left(\begin{array}{c|c}
\Phi^{(1)}_{11}&\Phi^{(1)}_{12}\\ \hline &\\[-.8pc]
\Phi^{(1)}_{21}&\Phi^{(1)}_{22}
\end{array}\right)$
which has a lower size.

\

\
 
\section{\shd The matrix functions $V_c$ and $W_c$ are left 
invertible in $H^\be$}
\setcounter{equation}{0}

\

In the last section we reduced the problem of finding a superoptimal solution
for $\Phi$ to the same problem for 
$\Phi^{(1)}\df\left(\begin{array}{c|c}
\Phi^{(1)}_{11}&\Phi^{(1)}_{12}\\\hline&\\[-1pc]\Phi^{(1)}_{21}&\Phi^{(1)}_{22}
\end{array}\right)$.
If we could continue this process, we would eventually reduce the problem to
the case $\min\{m_1,n_1\}=1$ and it would follow from Corollary \ref{t3.4}
that there is a unique superoptimal solution to the four block problem
for $\Phi$. The main problem now is to prove that the four block
operator $\G_{\Phi_1}$ has a maximizing vector. This is certainly the case
if $\|\G_{\Phi_1}\|_{\rm e}\le \|\G_\Phi\|_{\rm e}$.
To prove this inequality we use an idea of [PY2] based on the solution
of the so-called matricial corona problem for the matrix functions
$V_c$ and $W_c$. However in our case the solvability of this corona problem
is much harder than in [PY2] where $V_c,~W_c\in QC$.

In this section we shall prove that the matrix functions $V_c$ and $W_c$ are
left invertible in $H^\be$ (in other words the corona problem is solvable 
for them) which we shall use in the next section to prove that 
$\|\G_{\Phi_1}\|_{\rm e}\le \|\G_\Phi\|_{\rm e}$.

\begin{thm}\thd
\label{t5.1}
If $\|\G_\Phi\|_{\rm e}<\|\G_\Phi\|$, then the matrix functions $V_c$ 
and ${W_c}$ defined
in Section 3 are left invertible in $H^\be$.
\end{thm}

Clearly, it is sufficient to prove that $W_c$ is left invertible in $H^\be$,
which means that
there exists a matrix function $\O$ in $H^\be(M_{n_1,n_1-1})$ such that
$\O(\z)W_c(\z)=I$ for every $\z\in\dd$. To show the left invertibility
of $V_c$, it is sufficient to apply Theorem \ref{t5.1} to 
the transposed function $\Phi^t$
and use the equalities $\|\G_\Phi\|=\|\G_{\Phi^t}\|$ and
$\|\G_\Phi\|_{\rm e}=\|\G_{\Phi^t}\|_{\rm e}$, which follow immediately from 
the obvious identity
$$
\G_{\Phi^t}=J\G_\Phi^*J,
$$
where $J\r\df\bar z\bar\r$ for a vector function $\r$ in $L^2$.

Recall that $\overline w=\overline w_1\oplus\overline w_2$ 
is the first column of $W^*$.
Denote by $w_{1_r}$, $1\le r\le n_1$, the components of $w_1$.
We have $w_1=w_{(o)}w_{(i)}$, where
$w_{(o)}$ is a scalar outer function in $H^2$ and $w_{(i)}$ is an inner column 
function.

\begin{lem}\thd
\label{t5.2}
The vectorial Toeplitz operator $T_{\overline w_1}:H^2\to H^2(\C^{n_1})$ 
is left invertible.
\end{lem}

\pf First of all, $\Ker T_{\overline w_1}=\bO$. Indeed, assume that 
$\psi\in\Ker T_{\overline w_1}$. Then $\ov w_{1_r}\psi\in H^2_-$ for
$1\le r\le n_1$. Since ${\bf 1}$ is a greatest inner divisor of the components 
of $w_1$, it follows from Beurling's theorem that the functions
$$
\{\sum_{r=1}^{n_1}\kappa_rw_{1_r}:~\kappa_r\in H^2\}
$$
form a dense subset in $H^2$. Therefore we can approximate $\psi$ in the 
$L^1$-norm
by functions of the form 
$\sum_{r=1}^{n_1}\bar\kappa_r\bar w_{1_r}\psi$, each
of which belongs to 
$H^1_-\df\{\f\in L^1:~\hat \f(k)=0\mbox{ for } k\ge0\}$. Hence 
$\psi\in H^1_-$ and since $\psi\in H^2$, it follows that $\psi=\0$.

If $T_{\overline w_1}$ is not left invertible, there exists a sequence of scalar
functions $\{\f_j\}_{j\ge0}$ in $H^2$ such that $\|\f_j\|=1$ and 
$\f_n\to\0$ in the weak topology
and $\|T_{\overline w_1}\f_n\|\to0$. By Lemma \ref{t4.1} the operator
$T_{u_0}$ is onto and so there exists a sequence $\{\o_n\}_{n\ge0}$ of scalar 
functions in $(\Ker T_{u_0})^\perp$ such that $T_{u_0}\o_n=\f_n$.
Since $T_{u_0}$ is Fredholm, $\o_n\to\0$ weakly. Put 
$\r_n\df\o_nv\in H^2(\C^{n_1})\oplus L^2(\C^{n_2})$, 
where $v$ is the first 
column of $V$. Let $Q_0$ be an optimal solution of the four block problem
for $\Phi$. By (\ref{4.1}) we have
$$
\left(
\begin{array}{cc} \Phi_{11}-Q_0 & \Phi_{12} \\ 
                   \Phi_{21}& \Phi_{22} 
\end{array}
\right)\r_j=
W^*\left(\begin{array}{c}t_0u_0\o_j\\\0\\\vdots\\\0\end{array}\right)
=t_0u_0\o_j\overline w
=t_0(\f_j+\f_j^-)\overline w,
$$
for some functions $\f_j^-\in H^2_-$. It follows that
\bay
\|\G_{\Phi}\r_j\|^2&=&\left\|{\Bbb P}^-\left(
\begin{array}{cc} \Phi_{11}-Q_0 & \Phi_{12} \\ 
                   \Phi_{21}& \Phi_{22} 
\end{array}
\right)\r_j \right\|^2 = 
\|t_0{\Bbb P}^-(\f_j+\f_j^-)\overline{w}\|^2 \nonumber\\
&=& \|t_0(\f_j+\f_j^-) \overline{w}\|^2 
-\|\pp_+t_0(\f_j+\f_j^-) \overline{w}_1\|^2 \nonumber\\
 &= &\|t_0u_0\o_j\overline{w}\|^2-\|t_0\pp_+\f_j\overline{w}_1\|^2\nonumber\\
&=& \|t_0u_0\o_j\overline{v}\|^2 -  \|t_0\pp_+\f_j \overline{w}_1\|^2
=t_0^2(\|\r_j\|^2 -\|T_{\overline{w_1}}\f_j\|^2)\,, \nonumber
\ey
since $\|v(\z)\|_{\C^{n_1}}=\|w(\z)\|_{\C^{m_1}}$.

Taking into account that $\|T_{\overline{w_1}}\f_j\|\to0$ and $\r_j\to\0$ weakly,
we obtain 
$\|\G_\Phi\|_{\rm e}=t_0=\|\G_\Phi\|$ which contradicts the hypotheses of 
the lemma. $\bl$

The next step is to prove that the Toeplitz operator $T_{\ov w_{(i)}}$ is
left invertible, where $w_{(i)}$ is the inner part of $w_1$. We need the following 
well known facts. Let $\chi=\{\chi_j\}_{1\le j\le k}$ be a column function
in $H^\be(\C^k)$. Then it is left invertible in $H^\be$ (i.e. there exist functions
$\kappa_j$, $1\le j\le k$, such that $\sum_{j=1}^{k}\kappa_j(\z)\chi_j(\z)=1$
for all $\z\in\dd$) if and only if the Toeplitz operator $T_{\bar\chi}$ is left
invertible (see [Ar]). Note that by the Carleson corona theorem (see e.g., [Ni]) 
$\chi$ is left invertible if and only if $~\inf_{\z\in\dd}\|\chi(\z)\|_{\C^k}>0$.
This result was generalized in [SNF2] for the case of matrix (and even operator)
functions: let $\Xi$ be a matrix function in $H^\be$, then $\Xi$ is 
left invertible in $H^\be$
if and only if the Toeplitz operator $T_{\overline\Xi}$ is left invertible.

\begin{lem}\thd
\label{t5.3}
Under the hypotheses of Theorem \ref{t5.1} the Toeplitz operator
$T_{\ov w_{(i)}}$ is left invertible.
\end{lem}

\pf By Lemma \ref{t5.2}, $T_{\ov w_1}$ is left invertible. By Arveson's theorem
mentioned above
$w_1$ is left invertible in $H^\be$. We have $w_1=w_{(o)}w_{(i)}$, 
where $w_{(o)}$ is a scalar 
outer function in $H^\be$ and $w_{(i)}$ is an inner column function. Obviously,
it follows that $w_{(i)}$ is left invertible in $H^\be$. 
Again by Arveson's theorem
this implies that $T_{\ov w_{(i)}}$ is left invertible. $\bl$

We need the following result proved in [P].

\begin{thm}\thd 
\label{t5.4}
Let ${\bf W}$ be a unitary-valued matrix function of the form 
 ${\bf W}^t=
\left(\begin{array}{cc}{\bf w}& \overline W_c\end{array}\right)$, 
where ${\bf w}$ is a co-outer  inner column, and $W_c$ is a co-outer 
inner function. Then the Toeplitz operator $T_{{\bf W}^t}$ 
has trivial kernel and 
dense range, and the operators $H_{{\bf W}^t}^*H_{{\bf W}^t}$ and 
$H^*_{({\bf W}^t)^*}H_{({\bf W}^t)^*}$ are 
unitarily equivalent.
\end{thm}

The following result can easily be deduced from Theorem \ref{t5.4}

\begin{thm}\thd
\label{t5.5}
Let ${\bf W}$ be a matrix function satisfying the hypotheses of 
Theorem \ref{t5.4}.
Suppose that $\|H_{\overline{\bf w}}\|<1$. Then the Toeplitz operator 
$T_{{\bf W}^t}$ is invertible.
\end{thm}

\pf Clearly, 
$\|H_{({\bf W}^t)^*}\|=\|H_{{\bf w}^*}\|=\|H_{\ov{\bf w}^t}\|$. 
It is easy to see that
$\|H_{\ov{\bf w}^t} \|= \|H_{\ov{\bf w}}\|< 1$. 
By Theorem \ref{t5.4} the operators $H_{{\bf W}^t}^*H_{{\bf W}^t}$ and 
$H^*_{({\bf W}^t)^*}H_{({\bf W}^t)^*}$ 
are unitarily equivalent. Therefore $\|H_{{\bf W}^t}\| =\|H_{({\bf W}^t)^*}\|<1$.
Since ${\bf W}^t$ takes isometric values on $\T$, it is easy to see that 
$$
\|T_{{\bf W}^t} F\|_2^2+\|H_{{\bf W}^t}F\|_2^2=\|{\bf W}^tF\|_2^2=\|F\|_2^2
$$
for every vector function $F$. Consequently, $T_{{\bf W}^t}$ is left invertible 
if and
only if $\|H_{{\bf W}^t}\|<1$. It follows that both $T_{{\bf W}^t}$ 
and $T_{({\bf W}^t)^*}$ are left
invertible which means that $T_{{\bf W}^t}$ is invertible. $\bl$

{\bf Proof of Theorem \ref{t5.1}.} 
Put ${\bf w}\df w_{(i)}$ and let
${\bf W}^t=\left(\begin{array}{cc}w_{(i)}&\overline{W}_c\end{array}\right)$.
By Lemma \ref{t5.3}, $T_{\ov w_{(i)}}$
is left invertible. Since $w_{(i)}$ takes isometric values on $\T$,
we have as in the proof of Theorem \ref{t5.5}
$$
\|T_{\ov w_{(i)}}\o\|^2+\|H_{\ov w_{(i)}}\o\|^2=\|\o\|^2
$$
for every $\o\in H^2$. Hence $\|H_{\ov w_{(i)}}\|<1$ and so by
Theorem \ref{t5.5}
the operator $T_{{\bf W}^t}$ is invertible. Clearly, it follows that 
$T_{\ov W_c}$
is left invertible, since $T_{\ov W_c}$ can be interpreted as a restriction of
$T_{{\bf W}^t}$. 
Therefore by the Sz.-Nagy--Foias theorem mentioned above $W_c$ is
left invertible in $H^\be$. $\bl$

\

\

\section{\shd The essential norm of ${\bf\G_{\Phi^{(1)}}}$}
\setcounter{equation}{0} 

\

In Section 4 we reduced the proof of Theorem \ref{t2.3} to the fact that
 $\|\G_{\Phi_{(1)}}\|_{\rm e}\le\|\G_\Phi\|_{\rm e}$, where
the matrix function
$\Phi^{(1)}=\left(\begin{array}{c|c}
\Phi^{(1)}_{11}&\Phi^{(1)}_{12}\\\hline&\\[-1pc]\Phi^{(1)}_{21}&\Phi^{(1)}_{22}
\end{array}\right)$ is defined in (\ref{4.1}). In this section we are going
to use the facts that $V_c$ and $W_c$ are left invertible (see Section 5)
to prove this inequality which will complete the proof of Theorem \ref{t2.3}.

The idea behind the proof is 
the following. We use the fact that
$$
\|\G_{\Phi^{(1)}}\|_{\rm e}=\inf\{\limsup_{j}\|\G_{\Phi^{(1)}}\xi_j\|_2\},
$$
where the {\it infimum} is taken over all sequences $\{\xi_j\}$ in
$H^2(\C^{n_1-1})\oplus L^2(\C^{n_2})$ such that $\|\xi_j\|_2=1$ and
$\xi_j\to\0$ weakly. Given such a sequence $\{\xi_j\}$ we construct
another sequence $\{\r_j\}$ in $H^2(\C^{n_1})\oplus L^2(\C^{n_2})$
such that $\|\r_j\|_2=1$, $\r_j\to\0$ weakly and
$$
\limsup_j\|\G_\Phi\r_j\|_2\ge\limsup_{j}\|\G_{\Phi^{(1)}}\xi_j\|_2.
$$
To this end we are going to use a construction which is similar to the one used
in [PY2].

Let $W$ be the unitary-valued matrix function constructed in Section 3. Consider
the matrix $W^*$ which has the form
$$
W^*=\left(
\begin{array}{cc|c}
	\ov w_1 & W_c & F  \\
	\hline
	\ov w_2 & \0 & G  
\end{array}
\right).
$$
Let 
\begin{equation}
\label{6.1}
\b=\left(
\begin{array}{c|c}
	 W_c & F  \\
	\hline
	 \0 & G  
\end{array}
\right).
\end{equation}
To use a construction similar to the one given in [PY2], we 
have to find a left inverse of $\b$ of a special form. Recall that
we have proved in Section 5 that $W_c$ is left invertible in $H^\be$.
Let $W_c^{\rm li}$ be an $H^\be$ left inverse of $W_c$.

\begin{lem}\thd
\label{t6.1}
Let $\b$ be the matrix function defined by {\rm (\ref{6.1})}. Then the matrix function
$G$ is invertible in $L^\be$ and there exists a bounded 
left inverse of $\b$ of the form
\begin{equation}
\label{6.2}
B=\left(
\begin{array}{c|c}
	W_c^{\rm li} & X  \\
	\hline
	\0 & G^{-1}  
\end{array}
\right).
\end{equation}
\end{lem}

\pf Suppose that $G$ is not invertible in $L^\be$. Then there exists a sequence
$\{\xi_j\}_{j\ge0}$ in $L^2(\C^{m_2})$ such that $\|\xi_j\|_2=1$ and
$\|G\xi_j\|_2\to0$.

It is easy to see from Lemma \ref{t3.2} that
\begin{equation}
\label{6.3}
\|w_1(\z)\|^2_{\C^{m_1}}\ge\frac{\|\G_\Phi\|^2-\|\G_\Phi\|_{\rm e}^2}
{\|\G_\Phi\|^2}\df\d<1,~~~\z\in\T.
\end{equation}
Since the column $\ov w_1(\z)$ is orthogonal to the columns of $W_c(\z)$ a.e. on
$\T$, it follows that the matrix function 
$\left(\begin{array}{cc}\ov w_1 & W_c\end{array}\right)$ is invertible in $L^\be$.
Therefore there exists a bounded sequence $\{\eta_j\}$ in $L^2(\C^{m_1})$ such that
$$
\left(\begin{array}{cc}\ov w_1 & W_c\end{array}\right)\eta_j+F\xi_j=\0.
$$
Then
$$
W^*\left(\begin{array}{c}\eta_j\\\xi_j\end{array}\right)
=\left(\begin{array}{c}\0\\
\left(\begin{array}{cc}\ov w_2&\0\end{array}\right)\eta_j+G\xi_j\end{array}\right).
$$
It follows from (\ref{6.3}) that $\|w_2(\z)\|^2_{\C^{m_2}}\le1-\d$, $\z\in\T$.
Since $\|G\xi_j\|_2\to0$, we have for large values of $j$
$$
\left\|W^*\left(\begin{array}{c}\eta_j\\\xi_j\end{array}\right)\right\|_2
<\|\eta_j\|_2\le
\left\|\left(\begin{array}{c}\eta_j\\\xi_j\end{array}\right)\right\|_2,
$$
which contradicts the fact that $W$ is unitary-valued.

Let now $B$ be a matrix in the form (\ref{6.2}). Clearly $B\b=I$ if and only
if
$$
W_c^{\rm li}F+XG=\0.
$$ 
Since $G$ is invertible in $L^\be$, we can always find
a matrix function $X$ in $L^\be$ which satisfies this equality. $\bl$

{\bf Remark.} In the same way we can consider the submatrix $\a$ of the matrix
$\ov V$ constructed in Section 3,
$$
\a=\left(
\begin{array}{c|c}
	 V_c & \star  \\
	\hline
	 \0 & \star\end{array}\right)  
$$
and prove that $\a$ has a left inverse in the form
$$
A=\left(
\begin{array}{c|c}
	V_c^{\rm li} & \star  \\
	\hline
	\0 & \star  
\end{array}
\right),
$$
where $V_c^{\rm li}$ is an $H^\be$ left inverse of $V_c$.

To prove the main result of this section we need the following lemma
which in the case of Nehari's problem was proved in [PY2] (see Lemma 2.1 there).

\begin{lem}\thd
\label{t6.2}
Let $\eta$ be a vector function
in $H^2_-(\C^{m_1-1})\oplus L^2(\C^{m_2})$ and let $\chi$ be the scalar 
function in $H^2$ defined by
$$
\chi=-\pp_+w^tB^*\eta. 
$$
Then
$$
W^*\left(\begin{array}{cc}\chi\\\eta\end{array}\right)\in
H^2_-(\C^{m_1})\oplus L^2(\C^{m_2}).
$$
\end{lem}

\pf Since $W$ is unitary-valued, we have
$$
I=W^*W=\ov w w^t+\b\b^*
$$
and hence
$$
\b=\b(B\b)^*=\b\b^*B^*=(I-\ov w w^t)B^*.
$$
Therefore
$$
W^*\left(\begin{array}{cc}\chi\\\eta\end{array}\right)
=\left(\begin{array}{cc}\ov w&\b\end{array}\right)
\left(\begin{array}{cc}\chi\\\eta\end{array}\right)
=\ov w\chi+(I-\ov w w^t)B^*\eta=B^*\eta+\ov w(\chi+w^tB^*\eta).
$$
It is easy to see from (\ref{6.2}) that 
$B^*\eta\in H^2_-(\C^{m_1})\oplus L^2(\C^{m_2})$. Since 
 $w\in H^\be(\C^{m_1})\oplus L^\be(\C^{m_2})$, it follows that
$\ov w(\chi+w^tB^*\eta)=\ov w\pp_-w^tB^*\eta\in H^2_-$, which proves the result.
$\bl$

{\bf Remark.} It is easy to see that if $\{\eta_j\}$ is a sequence of functions in 
 $H^2_-(\C^{m_1-1})\oplus L^2(\C^{m_2})$ which converges weakly to $\0$, the
above construction produces a sequence of scalar functions $\{\chi_j\}$
in $H^2$, $\chi_j=-\pp_+(w^tB^*\eta_j)$, which also converges weakly to $\0$.

Now we are in a position to prove that 
$\|\G_{\Phi_1}\|_{\rm e}\le\|\G_\Phi\|_{\rm e}$, where the matrix function
$\Phi^{(1)}=\left(\begin{array}{c|c}
\Phi^{(1)}_{11}&\Phi^{(1)}_{12}\\\hline&\\[-1pc]\Phi^{(1)}_{21}&\Phi^{(1)}_{22}
\end{array}\right)$ is defined in (\ref{4.1}). 

\begin{thm}\thd
\label{t6.3}
Let $\Phi$ be a matrix function of the form {\rm (\ref{2.1})}
such that  $\|\G_\Phi\|_{\rm e}<\|\G_\Phi\|$.
Then $\|\G_{\Phi^{(1)}}\|_{\rm e}\le\|\G_\Phi\|_{\rm e}$.
\end{thm}

\pf Let $\{\xi_j\}$ be a sequence of functions in 
$H^2_-(\C^{n_1-1})\oplus L^2(\C^{n_2})$ such that $\|\xi_j\|_2=1$ and
$\xi_j\to\0$ weakly. Put $\eta_j=\G_{\Phi^{(1)}}\xi_j$. We are going to construct
a sequence of functions $\{\xi_j^{\#}\}$ in 
$H^2_-(\C^{n_1})\oplus L^2(\C^{n_2})$ such that 
$\frac{\xi_j^{\#}}{\|\xi_j^{\#}\|_2}\to\0$ weakly, and 
$\frac{\|\eta_j^{\#}\|_2}{\|\xi_j^{\#}\|_2}\ge\|\eta_j\|_2$,
where $\eta_j^{\#}\df\G_\Phi\xi_j^{\#}$.
As we have explained in the beginning of the section, this would imply
the desired inequality (put $\r_j=\frac{\xi_j^{\#}}{\|\xi_j^{\#}\|_2}$).

To this end we apply Lemma \ref{6.2} to the sequence $\{\eta_j\}$. We obtain 
a sequence of scalar $H^2$ functions $\{\chi_j\}$ such that
$$
W^*\left(\begin{array}{cc}\chi_j\\\eta_j\end{array}\right)\in
H^2_-(\C^{m_1})\oplus L^2(\C^{m_2}).
$$

Put
$$
\xi^{\#}_j=A^t\xi_j+q_jv,
$$
where $v$ is the first column of $V$ and $A$ is the left inverse of $\a$
described in the Remark after Lemma \ref{t6.1}. The scalar functions
$q_j$ will be chosen later. 

We have 
\begin{equation}
\label{6.4}
\left(\begin{array}{cc} t_0u_0 & \0 \\ \0 & \Phi^{(1)} \end{array}\right)
V^* \xi^{\#}_j =
\left(\begin{array}{c} t_0u_0 q_j +  t_0u_0 v^* A^t\xi_j\\ \Phi^{(1)}\xi_j 
\end{array}\right).
\end{equation}
Since  the Toeplitz operator $T_{u_0}$ is onto, we can pick $q_j$ as
a solution of the equation
$$
\pp_+(t_0u_0q_j+t_ou_ov^*A^t\xi_j)=\chi_j.
$$
Clearly, we may choose the $q_j$ so that $q_j\to\0$ weakly. Indeed, we may put
$$
q_j=t_0^{-1}(T_{u_0}|(\Ker T_{u_0})^\perp)^{-1}(t_0\chi_j-\pp_+u_0v^*A^t\xi_j).
$$
It follows that $\xi^{\#}_j\to\0$ weakly.

Let us show that the sequence $\{\xi_j^{\#}\}$ has the required 
properties. Since $\xi^{\#}_j\to\0$ weakly,
to prove that 
$\frac{\xi_j^{\#}}{\|\xi_j^{\#}\|_2}\to\0$ 
weakly, we have to estimate $\|\xi^{\#}_j\|_2$
from below. We have
\begin{equation}
\label{6.5}
\|\xi^{\#}_j\|_2^2=\|V^*\xi^{\#}_j\|_2^2=\|q_j+v^*A^t\xi_j\|^2_2+\|\xi_j\|^2_2\ge1.
\end{equation}
To complete the proof it remains to show that 
$\|\G_\Phi\xi_j^{\#}\|_2\ge\|\eta_j\|_2$.

Recall that $\eta_j=\G_{\Phi^{(1)}}\eta_j=\pp^-\Phi^{(1)}\xi_j$ and so
$\Phi^{(1)}\xi_j-\eta_j\in H^2(\C^{m_1-1})\oplus\{\0\}$. It is easy to see 
from the definition of $W$ (see Section 3) that
$$
W^*\left(\begin{array}{c}\0\\\Phi^{(1)}\xi_j-\eta_j\end{array}\right)
\in H^2(\C^{m_1})\oplus\{\0\}.
$$
It follows now from (\ref{6.4}) that
$$
\eta_j^{\#}=\pp^-\Phi\xi^{\#}_j=
\pp^-W^*\left(\begin{array}{cc}t_0u_0&\0\\\0&\Phi^{(1)}\end{array}\right)
V^*\xi_j^{\#}=\pp^-W^*\left(\begin{array}{c}\chi_j+\o_j\\\eta_j\end{array}\right),
$$
where
$$
\o_j\df\pp_-(t_0u_0q_j+t_0u_0v^*A^t\xi_j).
$$

Since the first column of $W^*$ is $\ov w_1\oplus\ov w_2$ and $w_1\in H^\be$,
it follows that
$$
W^*\left(\begin{array}{c}\o_j\\\0\end{array}\right)\in
\left(\begin{array}{c}H_-^2(\C^{m_1})\\L^2(\C^{m_2})\end{array}\right).
$$
We have chosen $\{\chi_j\}$ so that
$$
W^*\left(\begin{array}{c}\chi_j\\\eta_j\end{array}\right)\in
\left(\begin{array}{c}H_-^2(\C^{m_1})\\L^2(\C^{m_2})\end{array}\right).
$$
Therefore
$$
\pp^-W^*\left(\begin{array}{c}\chi_j+\o_j\\\eta_j\end{array}\right)
=W^*\left(\begin{array}{c}\chi_j+\o_j\\\eta_j\end{array}\right).
$$
Hence
$$
\|\eta_j^{\#}\|_2^2=\|\chi_j+\o_j\|_2^2+\|\eta_j\|^2_2=
t_0^2\|q_j+v^*A^t\xi_j\|^2_2+\|\G_{\Phi^{(1)}}\xi_j\|^2.
$$
Since $\|\G_{\Phi^{(1)}}\|\le t_0$, this together with (\ref{6.5})
yields 
$$
\frac{\|\eta_j^{\#}\|_2}{\|\xi_j^{\#}\|_2}\ge\frac{\|\eta_j\|_2}{\|\xi_j\|_2}=
\|\eta_j\|_2,
$$
which completes the proof. $\bl$

As we have already explained, Theorem \ref{t6.3} allows us to complete the 
proof of Theorem \ref{t2.3}.

{\bf Proof of Theorem \ref{t2.3}.} By Theorem \ref{t4.2} the four block for
$\Phi$ has a unique solution if so does the four block problem for $\Phi^{(1)}$ 
and the superoptimal singular values of the four block problem for
$\Phi^{(1)}$ are $t_1,t_2,\cdots,t_{d-1}$. By Theorem \ref{t6.3},
$\|\G_{\Phi^{(1)}}\|_{\rm e}\le\|\G_\Phi\|_{\rm e}$.
If $\G_{\Phi^{(1)}}=\0$, we certainly have uniqeness.  
Otherwise we can continue this process. 

Doing in this way we may stop
the process if we get on a certain stage the zero four block operator
or, otherwise, we eventually
reduce the problem to the case $d=1$. Uniqueness follows now from
Corollary \ref{t3.4}.

The fact that the singular values are constant on $\T$ follows immediately
from the facts that $V$ and $W$ are unitary-valued and $u_0$ is unimodular,
and from Lemma \ref{t3.3}. $\bl$

\

\

\section{\shd Thematic factorizations and indices of superoptimal singular
values}
\setcounter{equation}{0}

\

In this section we analyze the algorithm described in Section 3 and
obtain certain special factorizations of superoptimal symbols of 
four block operators satisfying the hypotheses of Theorem \ref{t2.3}.
Following [PY1] we shall call such factorizations 
{\it thematic}.

In Section 3 we have constructed matrix functions $V$ and $W$ associated with the
four block problem. By analogy with [PY1] we shall call matrix functions of the
form $V$ or $W$ {\it thematic functions}.

To state the result we may assume without loss of generality that $n_1\le m_1$
(otherwise we can take the transpose).

\begin{thm}\thd
\label{t7.1}
Let $\Phi$ be a superoptimal symbol of the four block operator
$\G$ which satisfies the hypotheses of Theorem
\ref{t2.3} and suppose that $~n_1\le m_1$.  Then
$\Phi$ admits the following factorization
\begin{equation}
\label{7.1}
\Phi=W^*_0 W^*_1W^*_2\cdots W^*_{d-1} D
V^*_{d-1}\cdots V^*_2 V^*_1 V^*_0
\end{equation} 
where
$$ D=\left(\begin{array}{cccc|c} t_0u_0 & \0 & \cdots & \0 &  \\ \0 & t_1u_1 &
\cdots &\0  &  \\
\vdots&\vdots&\ddots&\vdots & \\ \0 & \0 & \cdots & t_{d-1}u_{d-1} & \quad \0 
\quad \\
\0 & \0 & \cdots & \0      &     \\
\vdots&\vdots&\ddots&\vdots & \\  \0 & \0 & \cdots & \0  &   \\   \hline 
\multicolumn{4}{c|}{\0} & \star 
\end{array}\right),
$$
the $u_j$ are unimodular functions  such that the Toeplitz operator $T_{u_j}$
is Fredholm and $\ind T_{u_j}>0$, and the matrix functions $W_j$ and $V_j$
have the form
$$
V_j = 
\left(\begin{array}{cc}I_j & 0 \\ 0 & \breve V_j \end{array}\right)\, \qquad
W_j = 
\left(\begin{array}{cc}I_j & 0 \\ 0 & \breve W_j \end{array}\right),
$$ 
where $\breve V_j$, $\breve W_j$ are thematic matrix functions and $I_j$ is the identity
$j\times j$ matrix.
\end{thm}

It is easy to see that the successive application of the algorithm 
described in Section 3 gives us a desired factorization. 

{\bf Remark.} As in the case of Nehari's problem (see [PY1]) it is easy to see
that if a matrix function admits a factorization of the form (\ref{7.1}),
then it is the superoptimal symbol of the corresponding four block operator.

We can associate with
the factorization (\ref{7.1}) the {\it factorization indices} $k_j$
which are defined in the case $t_j\neq0$. 
We put $k_j=\ind T_{u_j}=\dim\Ker T_{u_j}$.

It was shown in [PY1] that even for Nehari's problem the indices depend 
on the choice of a thematic
factorization rather than on the function $\Phi$ itself. However, it was shown 
in [PY2]
that for Nehari's problem with compact Hankel operator 
the sum of the indices corresponding to equal
superoptimal singular values is an invariant (i.e. does not depend on the choice
of a factorization). 

The same turns out to be true for the four block problem too, and we shall 
prove
this later in Section 9.  Moreover, the sum of the indices 
corresponding to
equal superoptimal singular values admits a quite natural and simple 
geometric
interpretation.  To give this interpretation we have to introduce a new 
object
--- the so--called {\em superoptimal weight}.

\newpage

\section{\shd Superoptimal weight} 
\setcounter{equation}{0}

\

Let ${\cal W}\in L^\be(M_{n,n})$ be a {\em matrix weight}, i.e. a bounded
matrix-valued function on $\T$, whose values are nonnegative selfadjoint 
$n\times n$ matrices.

Given a four block operator 
$\G:H^2(\C^{n_1})\oplus L^2(\C^{n_2})\to 
H^2_-(\C^{m_1})\oplus L^2(\C^{m_2})$,
$n_1+n_2=n$, we call a weight $\W$ {\it admissible} 
if 
$$
\|\G f\|^2\le ({\cal W} f,f)\df\int_{\T}({\cal W}(\z)f(\z),f(\z))d{\bb m}(\z),
~~~  f\in H^2(\C^{n_1})\oplus L^2(\C^{n_2}).
$$

We  need the following result which we call {\em Generalized Nehari's
Theorem}.
\begin{thm}\thd
\label{t8.1}
Given a four block operator $\G_\Phi$ and an admissible weight $\W$
there exists a symbol $\Phi$ of $\G$ (i.e. an operator-valued function 
$\Phi$
such that $\G=\G_\Phi$) satisfying $\Phi^*\Phi\le \W$.
\end{thm}

If $\Phi$ is a symbol of $\G$ satisfying $\Phi^*\Phi\le \W$, we say
that $\Phi$ is {\it dominated by} the admissible weight $\W$.

In the case ${\cal W}\equiv cI$, $v\in\R_+$, this result was established in [FT],
and this is an analog of Nehari's theorem for four block operators.
In the general case the result follows from Theorem 1.1 of [TV], since
the four block operator $\G$ acting from the space 
$H^2(\C^{n_1}) \oplus L^2(\C^{n_2})$ endowed with the weighted norm 
$\|\cdot\|_{\W}$ 
to the space $H^2(\C^{m_1}) \oplus L^2(\C^{m_2})$ satisfies the hypothesis 
of the theorem. 

For the sake of completeness we deduce Theorem \ref{t8.1} from the analog
of Nehari's theorem mentioned above. 

{\bf Proof of theorem \ref{t8.1}.}
Define ${\cal W}_\e = {\cal W} + \e I$. Since ${\cal W}_\e \ge \e I$, 
it admits a factorization 
${\cal W}_\e=G_\e^*G_\e$, where $G_\e\in H^\infty(M_{n,n})$ 
is a matrix  function which is invertible in 
$H^\infty$ (see [R]). The weight $\W_\e$ is clearly admissible, so 
$$
\| \G f \| \le \| G_\e f\|\,,\qquad f \in H^2(\C^{n_1}) \oplus L^2(\C^{n_2})\,,
$$
which is equivalent to the fact that 
$$
\| \G G_\e^{-1} f \| \le \| f\|\,,\qquad f \in H^2(\C^{n_1}) \oplus L^2(\C^{n_2})\,,
$$
Since $G^{-1}_\e \in H^\infty(M_{n,n})$,
we can consider the operator $\G G_\e^{-1}$ as a four block operator. 
By the analog of Nehari's theorem it has a symbol $\Psi_\e$ such 
that $\|\Psi_\e\|_\infty\le1$. 
Then the function $\Phi_\e = \Psi G_\e$ is a symbol of $\G$ and 
$$
\Phi_\e^*\Phi_\e = G_\e^* \Psi_\e^*\Psi_\e G_\e \le G_\e^*G_\e = \W_\e =\W+\e I.
$$
It remains to chose a sequence $\{\e_j\}$ converging to 0 and such that
the sequence $\{\Phi_{\e_j}\}$ converges to a matrix function, say 
$\Phi\in L^\be$, in the $*$-weak topology. Clearly, $\Phi$ is a symbol of
$\G$ dominated by $\W$. $\bl$

{\bf Definition.} Let $\W$ be an admissible weight for the four block 
operator $G$. Consider the numbers
$$
s_j^\be({\cal W})\df\ess\sup_{\z\in\T}s_j({\cal W}(\z)),~~~0\le j\le d-1,
~~~d=\min\{m_1,n_1\}.
$$
The admissible weight $\W$ is called {\em superoptimal} if it 
lexicographically minimizes the numbers
$s_0^\be({\cal W}),\ s_1^\be({\cal W}), \cdots, s_{d-1}^\be(\W)$ 
among all admissible weights, i.e.,
$$
s_0^\be({\cal W})=\min\{s_0^\be({\cal V}):~{\cal V}\mbox{ is admissible}\},
$$
$$
s_1^\be({\cal W})=\min\{s_1^\be({\cal V}):~
{\cal V}\mbox{ is admissible},~s_0^\be({\cal V})\mbox{ is minimal possible}\},
~~~{\rm etc}.
$$

The following lemma shows that under the the hypotheses of Theorem \ref{t2.3} 
a superoptimal weight exists. However, a superoptimal weight is not 
unique in general. The lemma also shows  that  a superoptimal weight 
is nevertheless ``essentially'' unique for our purposes.  

Let $\l_a$, $a\in\R$, be the function on $\R$ defined by
$$
\l_a(t)\df\left\{
\begin{array}{ll} t,\quad & t\ge a \\[.5pc]
                  0,      & t < a
\end{array}\right..
$$

\begin{lem}\thd
\label{t8.2}
Let $\G$
be a four block operator satisfying the hypothesis of Theorem \ref{t3.2}.  
Let $\Phi$ be  
the superoptimal symbol of $\G$. Then 
\begin{enumerate}
\item $\Phi^*\Phi$ is a superoptimal weight for $\G$; 
\item If $\W$ and $\W'$ are two superoptimal weights, then $\l_a({\cal W}) = 
\l_a(\W')$ for any $a\ge t_{d-1}$. 
\end{enumerate}
\end{lem}

\pf Note that $\Phi$ is a symbol of $\G$ dominated by 
the weight $\Phi^*\Phi$. Suppose that $\Phi^*\Phi$ is not a 
superoptimal weight, i.e. that there exists an admissible weight $\W$ 
such that for some $j_0$, $0\le j_0\le d-1$,
$$
s_{j_0}^\be({\cal W})<s_{j_0}^\be(\Phi^*\Phi),~~~
s_j^\be({\cal W})=s_j^\be(\Phi^*\Phi), ~~~0\le j\le j_0.
$$
Let $\Psi$ be a symbol
of $\G$ dominated by the weight $\W$. Then 
$$
s_{j_0}^\be(\Psi)<s_{j_0}^\be(\Phi),~~~
s_j^\be(\Psi)=s_j^\be(\Phi), ~~~0\le j\le j_0,
$$ 
which contradicts the fact that $\Phi$ is the superoptimal symbol of $\G$.
Therefore $\Phi^*\Phi$ is a superoptimal weight. 

Let now $\W$ be a superoptimal weight, and let $\Psi$ be a symbol of $\G$ 
dominated by $\W$. Then $\Psi$ lexicographically 
minimizes
$(s_0^\infty(\Psi),\ s_1^\infty(\Psi), ... s_{d-1}^\infty(\Psi))$ and so
$\Psi$ coincides with the superoptimal symbol $\Phi$. So, for 
any superoptimal weight $\W$  the superoptimal symbol
$\Phi$ is the unique symbol of $\G$ dominated by $\W$. 
This means that $\Phi^*\Phi\le \W$ for any superoptimal weight $\W$. 
Together with the equalities $s_j^\be({\cal W})=
s_j^\infty(\Phi^*\Phi)=t_j^2$ this implies the second part of the lemma. 
$\bl$

Denote by $\L_a$  the function  defined by
\begin{equation}
\label{8.1}
\L_a(t)\df\left\{
\begin{array}{ll} t,\quad & t\ge a \\[.5pc]
                  a,      & t < a
\end{array}\right.
\end{equation}
The following fact is an easy consequence of Lemma \ref{t8.2}.

\begin{cor}\thd
\label{t8.3}
Let $\W$ and $\W'$ be two superoptimal weights. Then  $\L_a({\cal W}) = 
\L_a(\W')$ 
for any $a\ge t_{d-1}$.
\end{cor}

It is easy to see that if $a=t_{d-1}$ and $\W$ is a superoptimal weight, 
the weight $\L_a(\W)$ is the (unique) maximal superoptimal weight. 

\

\

\section{\shd Invariance of indices}
\setcounter{equation}{0}

\

The main result of this section shows that the sum of the indices of a
thematic factorization of a superoptimal symbol does not depend on the choice
of a factorization. To prove this fact we shall use the same construction
which was used in Section 6 to prove Theorem \ref{t6.3}.

Let 
$$
a_0>a_1>\cdots>a_l
$$
be all the distinct nonzero superoptimal singular values of a 
four block operator $\G$ which
satisfies the hypotheses of Theorem \ref{t2.3}. Let $\Phi$ be the superoptimal
symbol of $\G$ and let $k_j$ be the indices of a thematic factorization
of $\Phi$ of the form (\ref{7.1}). 
Consider the sum of the indices that correspond to equal superoptimal singular
values:
$$
\nu_r=\sum_{\{j:t_j=a_r\}}k_j,~~~0\le r\le l.
$$

The following theorem is the main result of the section.

\begin{thm}\thd
\label{t9.1}
The numbers $\nu_r$ do not depend on the choice of thematic factorization 
of $\Phi$.
\end{thm}

We are going to deduce Theorem \ref{t9.1} from Theorem \ref{t9.3} below, which
describes the numbers $\nu_r$ in terms of a superoptimal weight $\W$.

We say that a nonzero function $\xi\in H^2(\C^{n_1})\oplus L^2(\C^{n_2})$ 
is a {\it maximizing vector for an admissible weight} $\W$ if
$$
\|\G\xi\|^2 =({\cal W}\xi,\xi).
$$

\begin{lem}\thd
\label{t9.2}
Let $\G$ be a four block operator, $\W$  an admissible weight for $\G$, and 
$\Phi$  
a symbol of $\G$ dominated by $\W$. Let $\xi$ be a
maximizing vector for $\W$. 
Then $\Phi\xi\in H^2_-(\C^{m_1})\oplus L^2(\C^{m_2})$, i.e., 
$\G_\Phi\xi =\Phi\xi$. 
\end{lem}

Note that for ${\cal W}\equiv cI$, $c\in\R_+$, 
this was proved in Lemma \ref{t3.3}.

\pf We have
$$
({\cal W}\xi,\xi)= \|\G_\Phi\xi\|_2^2 =\|{\Bbb P}^-\Phi\xi\|_2^2 
\le \|\Phi\xi_2\|^2 \le (\W\xi,\xi)\,.
$$
It follows that $\|{\Bbb P}^-\Phi\xi\|_2^2\le \|\Phi\xi_2\|^2$, which implies the
result. $\bl$

Given an admissible weight $\W$  put 
$$
{\cal E}({\cal W})=\{\xi\in H^2(\C^{n_1})\oplus L^2(\C^{n_2}):~
\xi \mbox{ is maximizing for }\W \mbox{ or }\xi=\0\}.
$$
It is easy to see that $\xi\in\E(\W)$ if and only if
$\pp^+{\cal W}\xi-\G_\Phi^*\G_\Phi\xi=\0$, where $\pp^+$ is the orthogonal projection
onto $H^2(\C^{n_1})\oplus L^2(\C^{n_2})$. Therefore $\E(\W)$ is a closed linear
subspace. Recall that $\L_a(\W)$ does not depend on the choice of
superoptimal weight, where $\L_a$ is defined in (\ref{8.1}).

\begin{thm}\thd
\label{t9.3} 
Let $\G$ be a four block operator
satisfying the hypotheses of Theorem \ref{t2.3}, 
$\W$ a superoptimal admissible weight, and 
$\Phi$
the superoptimal symbol of $\G$. Consider a thematic factorization of $\Phi$
of the form {\rm (\ref{7.1})}. Let $k_j$ be the indices 
of the factorization. Then for $a\ge a_l$,
\begin{equation}
\label{9.1}
\sum_{\{j:t_j\ge a\}} k_j = \dim \E(\L_a(\W))\,.
\end{equation}
\end{thm}

Let us first deduce Theorem \ref{t9.1} from Theorem \ref{t9.3}.

{\bf Proof of theorem \ref{t9.1}.} It follows immediately from 
(\ref{9.1}) that 
$$
\nu_0=\dim\E(\L_{a_0}({\cal W})),~~~
\nu_j=\dim\E(\L_{a_j}({\cal W}))\ominus\E(\L_{a_{j-1}}({\cal W})),~~1\le j\le l,
$$
which proves the result. $\bl$

{\bf Proof of Theorem \ref{t9.3}.} It is easy to see that $\E(\L_a(\W))$
is constant on $(a_{j+1},a_j]$. So it is sufficient to prove that
for $0\le s \le l$
$$
\sum_{\{j:t_j\ge a_s\}} k_j = \dim \E({\cal W}_s)\,,
$$
where ${\cal W}_s=\L_{a_s}(\W)$. 

Let us prove the theorem by induction on $d$.

If $d=1$, factorization (\ref{7.1}) has the form
$$
\Phi=W^*_0DV^*_0\,,
$$
where 
$$
D=
\left(\begin{array}{c|c} t_0u_0 &     \\ 
\0 &   \\
\vdots& \0 \\ 
\0 &   \\
 \hline
{\0} & \star 
\end{array}\right),
$$
and $\|\star\|<t_0$ (otherwise the essential norm of $\G_\Phi$ would not be less than
$t_0$). Clearly, ${\cal W}_0(\z)\equiv a_0I$. If $\xi$ is a maximizing vector 
for ${\cal W}_0$, then it is easy to see that only the first entry of
$V^*_0\xi$ is nonzero. Therefore $\xi(\z)$ is pointwise orthogonal to all
columns of $V$ except for the first one. It follows that $\xi=\chi v$, where
$\chi$ is a scalar function in $L^2$. Using the fact that $v_1$ is a co-outer 
column function,
one can easily deduce that $\chi\in H^2$. It is easy to see that
$$
\Phi\xi=t_0u_0\chi\ov w.
$$
Since $\xi$ is a maximizing vector, it follows that
$\Phi\xi\in H^2_-(\C^{m_1})\oplus L^2(\C^{m_2})$. We can now use the fact that 
$w_1$ is a co-outer column function to deduce that $u_0\chi\in H^2_-$ which
means that $\chi\in\Ker T_{u_0}$. 

Conversely, it is easy to see that if
$\chi\in\Ker T_{u_0}$, then $\xi=\chi v$ is a maximizing vector, which proves 
that $\dim\E({\cal W}_0)=\dim\Ker T_{u_0}=k_0$.

Suppose now that the theorem is proved for $d-1$. We have
$$
\Phi=W^*\left(\begin{array}{cc}t_0u_0&\0\\\0&\Phi^{(1)}\end{array}\right)V^*,
$$
where $\Phi^{(1)}$ is the superoptimal symbol of $\G_{\Phi^{(1)}}$, 
$W\df W_0$, and
$V\df V_0$. The induction hypothesis implies that the theorem holds for 
$\G_{\Phi^{(1)}}$. 

Let $0\le s\le l$ and $a=a_s$. Suppose that $\W'$ is a superoptimal 
weight for $\G_{\Phi^{(1)}}$. By the induction hypothesis 
$$
\dim\E(\L_a({\cal W}'))=N\df\sum_{\{j\ge1:t_j\ge a\}}k_j,
$$
where the $k_j$ are the indices of the thematic factorization
(\ref{7.1}). By Lemma \ref{t9.2}, $\xi\in E(\L_a({\cal W}'))$ if and only if
$$
\Phi^{(1)}\xi \in H^2_-(\C^{m_1-1}) \oplus L^2(\C^{m_2})\qquad 
\mbox{and}\qquad \|\Phi^{(1)}\xi\|^2=(\L_a({\cal W}')\xi,\xi).
$$
Let $\xi_1,\xi_2,\cdots,\xi_N$ be a basis in $\E(\L_a({\cal W}'))$ and
let $\eta_\iota=\G_{\Phi^{(1)}}\xi_\i$. By Lemma \ref{t6.2} 
there exist scalar functions
$\chi_\i$, $1\le \i\le N$, such that
$$
 W^*\left(\begin{array}{c}\chi_\i \\ \eta_\i\end{array}\right) 
\in  H^2_-(\C^{m_1-1}) \oplus L^2(\C^{m_2}).
$$

As in the proof of Theorem \ref{t6.3} we define the functions $\xi^{\#}_\i$
as
$$
\xi^{\#}_\i=A^t\xi_\i+q_\i v,
$$
where $q_\i$ is a scalar functions in $H^2$ satisfying
$$
\pp_+(t_0u_0q_\i+t_0u_0v^*A^t\xi_\i)=\chi_\i
$$
(recall that the matrix function $A$ is defined in after Lemma \ref{t6.1}).
We have
$$
\eta^{\#}_\i\df {\Bbb P}^-\Phi \xi_\i^{\#} 
={\Bbb P}^-W^*
\left(\begin{array}{cc} t_0u_0 & \0 \\ \0 & \Phi^{(1)} \end{array}\right)
V^* \xi_\i^{\#}
= \pp^-W^*\left(\begin{array}{c}\chi_\i+\o_\i\\\eta_\i\end{array}\right),
$$
where as in the proof of Theorem \ref{t6.3}
$$
\o_\i=\pp_-(t_0u_0q_\i+t_0u_0v^*A^t\xi_\i).
$$
As we have explained in the proof of Theorem \ref{t6.3}
$$
\pp^-W^*\left(\begin{array}{c}\chi_\i+\o_\i\\\eta_\i\end{array}\right)=
W^*\left(\begin{array}{c}\chi_\i+\o_\i\\\eta_\i\end{array}\right)
$$
and so $\eta^{\#}_\i=\Phi\xi_\i^{\#}$.

Since the matrix function $W$ is unitary-valued, we have
\begin{eqnarray}
\|\Phi \xi^{\#}_\i\|^2 &=& \left\| W^* 
\left(
	\begin{array}{cc}
		t_0u_0 & \0  \\
		\0     & \Phi^{(1)} 
	\end{array}
\right)
\left(\begin{array}{c} q_\i + v^* A^t\xi_\i \\ \xi_\i 
\end{array} \right) \right\|^2 \nonumber \\[.5pc]
&=&t_0^2\|u_0(q_\i + v^* A^t\xi_\i) \|^2 +\|\Phi^{(1)}\xi_\i\|^2 \nonumber \\[.5pc] 
&=&t_0^2\|q_\i+v^* A^t\xi_\i\|^2+\left(\L_a({\cal W}')\xi_\i,\xi_\i\right)\nonumber 
\end{eqnarray}
(the last equality holds because $\xi_\i \in \E(\L_a(\W'))$, 
where ${\cal W}'=(\Phi^{(1)})^*\Phi^{(1)}$ is a superoptimal weight for 
$\G_{\Phi^{(1)}}$).

Consider the weight ${\cal V}$,
$$
{\cal V} = 
\left(
	\begin{array}{cc}
		t_0^2 & \0  \\
		\0     & \W' 
	\end{array}
\right).
$$
Bearing in mind that 
$$
V^*\xi_\i^{\#} = 
\left(\begin{array}{c} q_\i + v^* A^t\xi_\i \\ \xi_\i \end{array} \right), 
$$
we can continue the above chain of inequalities:
\begin{equation}
\label{9.2}
\|\Phi\xi^{\#}_\i\|^2 =(\L_a({\cal V}) V^* \xi^{\#}_\i, V^*\xi^{\#}_\i) =
(\L_a({\W}) \xi^{\#}_\i, \xi^{\#}_\i)
\end{equation}
(the last equality holds because $V$ is unitary-valued). Since
$\G_\Phi \xi^{\#}_\i ={\Bbb P}^- \Phi\xi^{\#}_\i = \Phi\xi^{\#}_\i$, 
it follows from 
(\ref{9.2}) that $\xi^{\#}_\i\in\E(\L_a(\W))$. 

We can add now another $k_0$ linear independent vectors of $\E(\L_a(\W))$.
Let $x_1,\cdots,x_{k_0}$ 
be a basis of $\Ker T_{u_0}$. Obviously,
$x_\i v\in\E(\L_a(\W))$. Let us show that the vectors 
$\xi^{\#}_1,\cdots,\xi^{\#}_N$, $x_1v,\cdots,x_{k_0}v$ are linearly 
independent. It is sufficient to prove that if $x\in\Ker T_{u_0}$ and
$xv+\sum_{\i=1}^Nc_\i\xi_\i^{\#}=\0$, then $x=\0$ and $c_\i=0$, $1\le\i\le N$. 
We have
\begin{equation}
\label{9.3}
V^*(xv+\sum_{\i=1}^Nc_\i\xi_\i^{\#})=\left(\begin{array}{c}x\\\0\end{array}\right)+
\sum_{\i=1}^Nc_\i\left(\begin{array}{c}v^*A^t\xi_\i+q_\i\\\xi_\i\end{array}\right)=\0.
\end{equation}
Since the $\xi_\i$ are linearly independent, it follows that
$c_\i=0$, $1\le\i\le N$, which in turn implies that $x=\0$.

This proves that 
$$
\sum_{\{j:t_j\ge a\}}k_j\le\dim\E(\L_a(\W)).
$$
Let us prove the opposite inequality. 

Denote by $\E_0$ the set of vectors in
$\E(\L_a(\W))$ of
the form $xv$ such that $x$ is a scalar function in $H^2$.
It is easy to see that $xv\in\E_0$ if and only if $x\in\Ker T_{u_0}$.
It remains to show that there exists at most $\sum_{\{j>0:t_j\ge a\}}k_j$
vectors $\breve\xi_\i$ that are linearly independent modulo $\E_0$. 
Let $\breve\eta_\i\df \G_\Phi \breve\xi_\i$. By Lemma \ref{t9.2},
$\breve\eta_\i =\Phi \breve\xi_\i$. Put
$$
V^*\breve\xi_\i = \left(\begin{array}{c} \gamma_j \\ \xi_\i \end{array}\right)\ 
,\qquad 
W^*\breve\eta_\i = \left(\begin{array}{c} \delta_\i \\ \eta_\i 
\end{array}\right)\ ,
$$
where $\gamma_\i$, $\delta_\i$ are scalar functions in $L^2$.  Since 
the  vectors $\breve\xi_\i$ are 
linearly independent modulo $\E_0$, the vectors $\xi_\i$ are linearly 
independent. To complete the proof, it is sufficient to show that $\xi_\i\in 
\E(\L_a(\W'))$. 

Since $\breve\eta_\i =\Phi \breve\xi_\i$, we have 
that $\eta_\i=\Phi^{(1)}\xi_\i$ and $\delta_\i=t_0u_0\gamma_\i$.  
It follows from the 
block structure of $V$ and $W$ that $\xi_\i\in 
H^2(\C^{n_1-1})\oplus L^2(\C^{n_2})$ and $\eta_\i\in 
H^2_-(\C^{m_1-1})\oplus L^2(\C^{m_2})$. 
So $\eta_\i = \Phi^{(1)}\xi_\i = \G_{\Phi^{(1)}}\xi_\i$. 

To show that $\xi_\i\in \E(\L_a(\W'))$, consider the following chain of 
equalities
\bay
\left(\L_a(\W)\breve\xi_\i, \breve\xi_\i\right)&=&
\left(\L_a(V{\cal \W} V^*)V\breve\xi_\i, V\breve\xi_\i\right)\nonumber\\[.5pc]
&=& \left(
	\left(
	\begin{array}{cc}
		t_0^2 & \0  \\
		\0     & \L_a(\W')
	\end{array}
	\right)
	\left(
	\begin{array}{c}
		\gamma_\i  \\
		 \xi_\i
	\end{array}
	\right),
	\left(
	\begin{array}{c}
		\gamma_\i  \\
		 \xi_r
	\end{array}
	\right)
\right)\nonumber
\\[.5pc]
&=&
t_0^2\|\gamma_\i\|^2+\left(\L_a({\cal W}')\xi_\i, \xi_\i\right).\nonumber
\ey
On the other hand 
$$
\left(\L_a({\cal W})\breve\xi_\i,\breve\xi_\i\right)= \|\Phi \breve\xi_\i\| =
\|\breve\eta_\i\| = \|\eta_\i\|^2+\|\delta_\i\|^2= \|\Phi^{(1)}\xi_\i\|^2 + 
t_0^2\|\gamma_\i\|^2.
$$
Therefore
$\left(\L_a({\cal W}')\xi_\i, \xi_\i\right) = \|\Phi^{(1)}\xi_\i\|^2$,
which implies
$\xi_\i\in\L_a(\W')$. $\bl$

\

\

\section{\shd Singular values of ${\bf \G_\Phi}$ and 
superoptimal singular values}
\setcounter{equation}{0}

\

Let $\G$ be a four block operator satisfying the hypotheses of Theorem
\ref{t2.3}. Denote by $\Phi$ its unique superoptimal symbol and
consider a thematic factorization of $\Phi$ of the form (\ref{7.1}).
Let $\{t_j\}$ be the superoptimal singular values and
$\{k_j\}$ the indices of the factorization. Consider the extended 
$t$-sequence for $\G$:
$$
t_0,t_0,\cdots,t_1,t_1,\cdots,t_1,\cdots
$$
in which $t_j$ is repeated $k_j$ times. We denote the terms of the extended 
sequence by $t'_0,t'_1,t'_2,\cdots$. 
Although the indices $k_j$ depend on the choice of thematic factorization,
it follows from Theorem \ref{9.1} that the extended $t$-sequence is uniquely
determined by $\G$.

In [PY2] it was shown in the case of Nehari's problem with $\Phi\in H^\be+C$ that
$t'_j\le s_j(H_\Phi)$, $0\le j\le k-1$. In this section we are going to prove
the same inequality in the case of the four block problem under the hypotheses of
Theorem \ref{t2.3}. Moreover, we prove in this section a stronger result
which is also new in the case of Nehari's problem with an $H^\be+C$ symbol.
To prove the results we use in this section the same machinery as we used in 
Section 6.

Let $\G$ be a four block operator that satisfies the hypotheses of Theorem 
\ref{t2.3}. Then (see Section 3)
$$
\Phi=W^*\left(
\begin{array}{cc}t_0u_0&\0\\\0&\Phi^{(1)}\end{array}\right)V^*,
$$
where the unitary-valued matrix functions $V$ and $W$ are defined in
Section 3. The following inequality is the main result of the section.

\begin{thm}\thd
\label{t10.1}
Let $\G$ be a four block operator such that
$\|\G_\Phi\|_{\rm e}<\|\G_\Phi\|$
and let $\Phi$ be its superoptimal symbol. Then
$$
s_j(\G_{\Phi^{(1)}})\le s_{j+k_0}(\G_\Phi),~~~j\in\Z_+.
$$
\end{thm}

Recall that $k_0=\dim\Ker T_{u_0}$.

Let us first derive from Theorem \ref{t10.1} the desired inequality between
the singular values of $\G_\Phi$ and the superoptimal singular values.

\begin{thm}\thd
\label{t10.2}
Under the hypotheses of Theorem \ref{t2.3}
$$
t'_j\le s_j(\G_\Phi),~~~j\ge0.
$$
\end{thm}

{\bf Proof of Theorem \ref{t10.2}.}  Let $x\in\Ker T_{u_0}$. Clearly,
$$
\G_\Phi xv=W^*\left(\begin{array}{cccc}t_0u_0x&\0&\cdots&\0\end{array}\right)^t
=t_0u_0x\ov w\in H^2_-(\C^{m_1})\oplus L^2(\C^{m_2}).
$$
It follows that $\|\G_\Phi xv\|_2=t_0\|xv\|_2$ which proves that
\begin{equation}
\label{10.1}
s_j(\G_\Phi)=t_0,~~~0\le j\le k_0-1.
\end{equation}

We can now proceed by induction on $d$. Clearly, the result holds for $d=1$.
It is also obvious that if the theorem holds for $\G_{\Phi^{(1)}}$, then
by Theorem \ref{t10.1}
$$
t'_j\le s_j(\G_\Phi),~~~k_0\le j\le d-1,
$$
which together with (\ref{10.1}) proves the theorem. $\bl$

{\bf Proof of Theorem \ref{t10.1}.} Clearly, it is sufficient to prove the 
following fact. Let $\cL$ be a subspace of 
$H^2_-(\C^{n_1-1})\oplus L^2(\C^{n_2})$ such that 
$\|\G_{\Phi^{(1)}}\xi\|_2\ge s\|\xi\|_2$, for every $\xi\in\cL$, 
where $0<s\le t_0$, 
then there exists a subspace $\M$ of
$H^2_-(\C^{n_1})\oplus L^2(\C^{n_2})$ such that $\dim\M\ge\dim\cL+k_0$ and
$\|\G_\Phi\r\|_2\ge s\|\r\|_2$ for every $\r\in\M$.

Let $\xi_\i$, $1\le\i\le N$, be a basis in $\cL$. Put 
$\eta_\i=\G_{\Phi^{(1)}}\xi_\i$. By Lemma \ref{6.2} there exist scalar functions
$\chi_\i$ in $H^2$ such that 
$W^*\left(\begin{array}{c}\chi_\i\\\eta_\i\end{array}\right)\in 
H^2_-(\C^{m_1})\oplus L^2(\C^{m_2})$. We define the functions
$\xi^{\#}_\i\in H^2_-(\C^{n_1})\oplus L^2(\C^{n_2})$ by
$$
\xi_\i^{\#}=A^t\xi_\i+q_\i v,
$$
where $q_\i$ is a scalar function in $H^2$ satisfying
$$
\pp_+(t_0u_0q_\i+t_0u_0v^*A^t\xi_\i)=\chi_\i
$$
(see the proof of Theorem \ref{t6.3}). 

We can now define $\M$ by
$$ 
\M=\spn\{\xi_\i^{\#}+xv:~1\le\i\le N,~x\in\Ker T_{u_0}\}.
$$

Let us show that $\dim\M=N+k_0$. Since $\dim\Ker T_{u_0}=k_0$, it is sufficient 
to prove that if $xv+\sum_{\i=1}^Nc_\i\xi^{\#}_\i=\0$, then $x=\0$ and $c_\i=0$,
$1\le\i\le N$. This follows immediately from (\ref{9.3}).

To complete the proof it remains to show that $\|\G_\Phi \r\|_2\ge s\|\r\|_2$
for  $\r=xv+\sum_{\i=1}^rc_j\xi_\i^{\#}$. 
Let $\xi=\sum_{\i=1}^rc_\i\xi_\i$, 
$\eta=\G_\Phi^{(1)}\xi$,
$q=\sum_{\i=1}^rc_\i q_\i$, and 
 $\xi^{\#}_\i=\sum_{\i=1}^rc_\i\xi^{\#}_\i$. 

We have
\begin{eqnarray*}
W^*\left(\begin{array}{cc}t_0u_0&\0\\\0&\Phi^{(1)}\end{array}\right)\r&=&
W^*\left(\begin{array}{cc}t_0u_0&\0\\\0&\Phi^{(1)}\end{array}\right)
(xv+qv+A^t\xi)\\[.5pc]
&=&W^*\left(\begin{array}{cc}t_0u_0&\0\\\0&\Phi^{(1)}\end{array}\right)
\left(\begin{array}{c}x+q+v^*A^t\xi\\\xi\end{array}\right)\\[.5pc]
&=&W^*\left(\begin{array}{c}t_0u_0x+t_0u_0q+t_0u_0v^*A^t\xi\\\Phi^{(1)}\xi
\end{array}\right).
\end{eqnarray*}
It follows (see the proof of Theorem \ref{t6.3}) that
$$
\G_\Phi\r=
W^*\left(\begin{array}{c}t_0u_0x+t_0u_0q+t_0u_0v^*A^t\xi\\\G_{\Phi^{(1)}}\xi
\end{array}\right).
$$
Therefore
$$
\|\G_\Phi\r\|^2_2=|t_0|^2\|x+q+v^*A^t\xi\|^2_2+\|\eta\|^2_2.
$$
We have
$$
\|\r\|_2^2=\|V^*\r\|^2_2=\|x+q+v^*A^t\xi\|_2^2+\|\xi\|^2_2.
$$
Since $s\le t_0$ and $\|\eta\|_2\ge s\|\xi\|_2$, 
it follows that $\|\G_\Phi\r\|^2_2\ge s^2\|\r\|^2_2$.
$\bl$

Theorem \ref{t10.1} certainly applies to the case of Nehari's
problem. Recall that a matrix function 
$\Phi$ is called {\it very badly approximable} (see [PY]) if the zero function
is a superoptimal approximant of $\Phi$. 

Recall that under the condition $\|H_\Phi\|_{\rm e}<\|H_\Phi\|$ the function
$\Phi$ admits a factorization
\begin{equation}
\label{10.2}
\Phi=W^*\left(
\begin{array}{cc}t_0u_0&\0\\\0&\Phi^{(1)}\end{array}\right)V^*,
\end{equation}
where $V$ and $W$ are unitary matrix functions of the form
$$
V=\left(\begin{array}{cc}v&\ov V_c\end{array}\right),~~~~
W=\left(\begin{array}{cc}w&\ov W_c\end{array}\right)^t,
$$
and $u_0$ is a unimodular function such that $k_0\df\dim\Ker T_{u_0}>0$.

The following result is certainly
a partial case of Theorem \ref{t10.1}.

\begin{thm}\thd
\label{t10.3}
Let $\Phi$ be a very badly approximable matrix function on $\T$ such
that $\|H_\Phi\|_{\rm e}$ is less that the smallest nonzero
superoptimal singular value of Nehari's problem. 
Then
$$
s_j(H_{\Phi^{(1)}})\le s_{j+k_0}(H_\Phi),~~~j\in\Z_+,
$$
where $\Phi^{(1)}$ and $k_0$ are given by the factorization {\rm (\ref{10.2})}.
\end{thm}

\

\end{document}

%% file: lastart.tex
\renewcommand{\a}{\alpha}
\renewcommand{\b}{\beta}

\renewcommand{\d}{\delta}
\newcommand{\e}{\varepsilon}
\newcommand{\z}{\zeta}

\newcommand{\vt}{\vartheta}
\renewcommand{\l}{\lambda}
\renewcommand{\r}{\rho}

\renewcommand{\t}{\tau}

\newcommand{\f}{\varphi}
\renewcommand{\o}{\omega}
\newcommand{\w}{\omega}
\newcommand{\G}{\Gamma}

\renewcommand{\L}{\Lambda}

\renewcommand{\O}{\Omega}

\newcommand{\E}{{\cal E}}

\newcommand{\cL}{{\cal L}}
\newcommand{\M}{{\cal M}}

\newcommand{\C}{{\Bbb C}}
\newcommand{\T}{{\Bbb T}}
\newcommand{\pp}{{\Bbb P}}
\newcommand{\dd}{{\Bbb D}}
\newcommand{\R}{{\Bbb R}}
\newcommand{\Z}{{\Bbb Z}}
\newcommand{\0}{{\Bbb O}}
\newcommand{\bO}{\{ {\Bbb O} \}}

\newcommand{\bb}[1]{\mbox{{\boldmath$#1$\unboldmath}}}

\newcommand{\df}{\stackrel{\rm{def}}{=}}
\newcommand{\dist}{\operatorname{dist}}
\newcommand{\Ker}{\operatorname{Ker}}

\newcommand{\spn}{\operatorname{span}}

\newcommand{\beq}{\begin{equation}}
\newcommand{\eeq}{\end{equation}}
\newcommand{\bay}{\begin{eqnarray}}
\newcommand{\ey}{\end{eqnarray}}

\newcommand{\be}{\infty}

\newcommand{\bl}{\blacksquare}
\newcommand{\ess}{\operatorname{ess}}
\newcommand{\ind}{\operatorname{ind}}

\newcommand{\Range}{\operatorname{Range}}

\newtheorem{thm}{\hspace{\parindent}Theorem}[section]

\newtheorem{cor}[thm]{\hspace{\parindent}Corollary}
\newtheorem{lem}[thm]{\hspace{\parindent}Lemma}